\documentclass[10pt, a4paper]{article}
\usepackage{amsmath,amssymb,amsfonts,cite,amsthm}
\usepackage{bm}
\usepackage{indentfirst}
\usepackage{pdflscape}
\usepackage{graphicx}
\usepackage{tabularx,array,hhline,booktabs}
\usepackage{enumerate}
\usepackage{setspace}
\usepackage{xcolor, color, colortbl}
\usepackage{multirow}
\usepackage{dcolumn}
\usepackage{textcomp}
\newcolumntype{d}[1]{D{.}{.}{#1}}
\newcolumntype{M}[1]{>{\centering\arraybackslash}m{#1}}
\newcolumntype{N}{@{}m{0pt}@{}}
\usepackage[british]{babel}
\usepackage{authblk}
\usepackage{hyperref}
\usepackage{algorithm}
\usepackage{algorithmicx}
\usepackage{algpseudocode}
\usepackage{siunitx}
\usepackage{enumitem}
\usepackage{caption}

\numberwithin{equation}{section}
\numberwithin{figure}{section}
\numberwithin{table}{section}

\newtheorem{remark}{Remark}
\numberwithin{remark}{section}

\newcolumntype{d}[1]{D{.}{.}{#1}}
\newcolumntype{M}[1]{>{\centering\arraybackslash}m{#1}}
\newcolumntype{N}{@{}m{0pt}@{}}
\DeclareMathOperator*{\argmin}{\arg\min}

\definecolor{gray1}{gray}{0.95}
\definecolor{gray2}{gray}{0.85}
\definecolor{gray3}{gray}{0.75}
\definecolor{gray4}{gray}{0.65}

\newcommand\kwd[1]{\par{\textbf{\textit{Keyword:}}~}{#1}}

\title{Laplacian regularized eikonal equation with Soner boundary condition on polyhedral meshes}
\author[1]{Jooyoung Hahn\footnote{jooyoung.hahn@stuba.sk}}
\author[1]{Karol Mikula\footnote{karol.mikula@stuba.sk}}
\author[1]{Peter Frolkovi\v{c}\footnote{peter.frolkovic@stuba.sk}}
\affil[1]{Faculty of Civil Engineering, Slovak University of Technology, Department of Mathematics and Descriptive Geometry, Radlinsk{\'e}ho 11, 810 05 Bratislava, Slovak Republic}
\date{} %% if you don't need date to appear
\begin{document}
	\maketitle
	\begin{abstract}
		In this paper, we propose a numerical algorithm based on a cell-centered finite volume method to compute a distance from given objects on a three-dimensional computational domain discretized by polyhedral cells. Inspired by the vanishing viscosity method, a Laplacian regularized eikonal equation is solved and the Soner boundary condition is applied to the boundary of the domain to avoid a non-viscosity solution. As the regularization parameter depending on a characteristic length of the discretized domain is reduced, a corresponding numerical solution is calculated. A convergence to the viscosity solution is verified numerically as the characteristic length becomes smaller and the regularization parameter accordingly becomes smaller. From the numerical experiments, the second experimental order of convergence in the $L^1$ norm error is confirmed for smooth solutions. Compared to solve a time-dependent form of eikonal equation, the Laplacian regularized eikonal equation has the advantage of reducing computational cost dramatically when a more significant number of cells is used or a region of interest is far away from the given objects. Moreover, the implementation of parallel computing using domain decomposition with $1$-ring face neighborhood structure can be done straightforwardly by a standard cell-centered finite volume code.
	\end{abstract}
	\kwd{Vanishing viscosity method, Eikonal equation, Soner boundary condition, Laplacian regularizer, Cell-centered finite volume method, Polyhedral meshes}
	
\section{Introduction}
The viscosity solution of an eikonal equation is used in various applications from pure geometrical analysis to complicated problems mentioned in~\cite{ref:S99book,ref:OF00}. In the premixed turbulent combustion with thin flame fronts~\cite{ref:P00}, a distance from the thin flame modeled by a surface is used to design the flame-wall interaction and quenching~\cite{ref:SLSE17} or the end-gas autoignition for knock prediction~\cite{ref:M16}. A distance from a computational boundary, so-called wall distance, is a crucial feature in turbulence modeling methods~\cite{ref:BL78,ref:BB91,ref:SA94,ref:FS02,ref:T03}. It is also useful to obtain the medial axis transformation~\cite{ref:XT10,ref:XT11} of a given domain, which is crucial to automated mesh generation~\cite{ref:PAS95,ref:QRPG04}. In cardiac electrophysiology~\cite{ref:CFG93,ref:K91,ref:THP02}, a properly modeled eikonal equation approximates a propagation of excitation wavefront by the time to excite all points in the myocardium. In geophysics, a propagation of seismic waves is described by an eikonal equation in the high frequency regions~\cite{ref:RS05}.

In order to make more realistic simulation of the mentioned applications, it is necessary to use three-dimensional (3D) discretized domain in a fine scale to capture detailed phenomena. On such a domain, a parallel computing using domain decomposition is inevitable because of significantly high consumption of the memory. Moreover, computational domains of the industrial problems described by a complicated boundary shape are commonly discretized by polyhedral cells; see more advantages to use polyhedral cells~\cite{ref:P04}. Therefore, the target we would like to achieve here is to compute a distance function from given objects on polyhedral meshes by a parallel computing using domain decomposition with the simplest structure of overlapping domains, that is, $1$-ring face neighbor structure~\cite{ref:HMFB17}. For usability of the developed algorithm, it should be possible to make a straightforward implementation in a standard code of cell-centered finite volume method (FVM). 

The most well-known algorithm to efficiently solve an eikonal equation is usually considered to be the fast marching method (FMM)~\cite{ref:S96,ref:BS98,ref:KS98}. The fast computation is obtained by keeping a heap data structure to handle active nodes on a propagating front as candidates for updating the values. However, for typical parallel computing using domain decomposition, the heap structure is difficult to be maintained efficiently in parallel computation. An alternative approach is the fast sweeping method (FSM)~\cite{ref:Z05,ref:Z07,ref:QZZ07a,ref:TCOZ03} by updating necessary values with a Gauss-Seidel type iterations and it achieves better computational speed in a simple computational domain because a sorting is not used; see detailed computational study of FMM and FSM in~\cite{ref:HT05,ref:GK06}. In the fast iterative method (FIM)~\cite{ref:JW08,ref:FJPKW11,ref:FKW13}, a fine-grained parallel algorithm to solve an eikonal equation is presented on regular square, triangular, and tetrahedron meshes. A particular assumption to use FIM and FMM on triangular or tetrahedral meshes is that a shape of cell is restricted to an acute triangle or tetrahedron. For obtuse shapes, a smart division is necessary to make all cells as acute shapes but it is not clear how efficiently it can be divided in polyhedral meshes in a situation of moving mesh or remeshing that commonly happens in combustion simulation. 

In this paper, we numerically find a viscosity solution of an eikonal equation:
\begin{equation}\label{eq:EkEq_org}
	\begin{alignedat}{2}
		|\nabla u(\mathbf{x})| &= 1, \quad &&\mathbf{x} \in \Omega\setminus\Gamma, \\
		u(\mathbf{x}) &= 0, \quad &&\mathbf{x} \in \Gamma,
	\end{alignedat}
\end{equation}
where a computational domain~$\Omega$ is either convex or non-convex and~$\Gamma$ indicates fixed locations represented by a collection of curves or surfaces or a part of the boundary of the computational domain. The viscosity solution of~\eqref{eq:EkEq_org} defined in~~\cite{ref:CEL84} is the distance function from~$\Gamma$ on the domain~$\Omega$. A noticeable necessary condition of being the viscosity solution of~\eqref{eq:EkEq_org} is an inequality condition on the boundary of the domain:
\begin{equation}\label{eq:SonerBC}
	\bm{\nu}(\mathbf{x}) \cdot \nabla u(\mathbf{x}) \geq 0, \quad \mathbf{x} \in \partial \Omega \setminus\Gamma,
\end{equation}
where $\bm{\nu}$ is the outward normal to the boundary of the domain. The above inequality is proved in the \textit{Remark} after Proposition II.1 in \cite{ref:CDL90}. It is so-called the Soner boundary condition~\cite{ref:DES11} or the state constraint condition in optimal control problems~\cite{ref:CDL90,ref:S86}.
The condition is applied on obstacle boundaries~\cite{ref:FT09} and it restricts the discrete set of admissible control on all points in a domain in order to avoid an incorrect search direction. A general shape of obstacle embedded in a discretized domain is considered in~\cite{ref:GK06}. The eikonal equation~\eqref{eq:EkEq_org} and the Soner boundary condition~\eqref{eq:SonerBC} are discretized by a monotone finite difference scheme in~\cite{ref:DES11} when a set $\Gamma$ is a collection of finite discrete points and the error bound of the scheme is derived to the order of the square of cell size on a regular rectangular mesh. The obstacle~\cite{ref:FT09} can be understood as a hole in a domain~\cite{ref:HMFB21} and the necessity of using the Soner boundary condition and its geometrical interpretation are explained in~\cite{ref:HMFB21} by numerical examples.

A time-relaxed formulation of~\eqref{eq:EkEq_org} with the Soner boundary condition~\eqref{eq:SonerBC} is presented to compute a signed distance function when a shape of $\Gamma$ is a closed, bounded, orientable, and connected surface $\Gamma$ in a general computational domain $\Omega \subset \mathbb{R}^3$~\cite{ref:HMFB21}:
\begin{equation}\label{eq:EkEq_Time}
	\begin{alignedat}{2}
		\frac{\partial}{\partial t} \phi(\mathbf{x},t) \pm |\nabla \phi(\mathbf{x},t)| &= \pm 1 \quad && (\mathbf{x},t) \in \Omega^{\pm}\times(0,T], \\
		\phi(\mathbf{x},t) &= 0 \quad && (\mathbf{x},t) \in \Gamma\times[0,T], \\
		\bm{\nu}(\mathbf{x}) \cdot \nabla \phi(\mathbf{x},t) & \geq 0 \quad && (\mathbf{x},t) \in (\partial \Omega \setminus \Gamma) \times(0,T],
	\end{alignedat}
\end{equation}
where $\phi(\mathbf{x},0) > 0$ on $\Omega^{+}$ and $\phi(\mathbf{x},0) < 0$ on $\Omega^{-}$ are outside and inside the closed surface, respectively. The Soner boundary condition is essential to avoid a non-viscosity solution, especially on a non-convex domain. The distance information from $\Gamma$ is propagated into the rest of domain $\Omega \setminus \Gamma$ along the normal direction to $\Gamma$ over the time. The steady state solution eventually becomes a signed distance function from $\Gamma$. In the case of computing a wall distance function, that is, $\Gamma = \partial \Omega$, a transport form of eikonal equation~\eqref{eq:EkEq_org} is presented in~\cite{ref:XT10} and the algorithm is implemented by a standard FVM code with the first order upwind scheme. Even if the time relaxation in~\cite{ref:XT10,ref:HMFB21} with a proper choice of time step brings a robustness of the algorithm, a main disadvantage of using~\eqref{eq:EkEq_Time} is a large amount of computational cost when a region of interest is located far away from $\Gamma$.

Inspired by the vanishing viscosity method~\cite{ref:CDL90}, the equation we numerically solve is combined with a Laplacian regularizer and the solution $u_\epsilon$ is an approximation of the viscosity solution of~\eqref{eq:EkEq_org}:
\begin{equation}\label{eq:EkEq_Lap}
	\begin{alignedat}{2}
		-\epsilon \triangle u_\epsilon(\mathbf{x}) + |\nabla u_\epsilon(\mathbf{x})| &= 1 \quad &&\mathbf{x} \in \Omega\setminus\Gamma, \\
		u_\epsilon(\mathbf{x}) &= 0 \quad && \mathbf{x} \in \Gamma, \\
		\bm{\nu}(\mathbf{x}) \cdot \nabla u_\epsilon(\mathbf{x}) & \geq 0 \quad && \mathbf{x} \in \partial \Omega \setminus\Gamma,
	\end{alignedat}
\end{equation}
where $\epsilon > 0$ is the regularization parameter. Compared to solve~\eqref{eq:EkEq_Time}, a clear advantage of solving the above equation is to improve computational cost because of an infinite propagation speed caused by the Laplacian regularization term. In order to numerically solve~\eqref{eq:EkEq_Lap}, two difficulties should be resolved: the first is how to deal with the nonlinear term and the second is how to choose a regularization parameter.  In~\cite{ref:TRSBB05,ref:T11}, the same Laplacian regularizer is used for computing a wall distance function, that is, $\Gamma = \partial \Omega$. The non-linearity in~\eqref{eq:EkEq_Lap} is resolved by using $|\nabla u_{\epsilon}|^2$ and its linearization. The choice of the regularization parameter depends on an approximated distance from $\Gamma$, which makes more inaccurate results on the far field. In~\cite{ref:BF15,ref:CGP15,ref:BF19,ref:EIN21}, the non-linearity in~\eqref{eq:EkEq_org} is managed by an energy minimization with the constraint $\mathbf{p} = \nabla u$ and then a penalty method or augmented Lagrangian method are used to approximate a viscosity solution of~\eqref{eq:EkEq_org} for the cases of $\Gamma = \partial \Omega$. Throughout this paper, we discuss the details of two mentioned difficulties of solving~\eqref{eq:EkEq_Lap} in order to obtain a meaningful convergence order numerically.

The rest of paper is presented as follows. In Section~\ref{sec:prop_alg}, we explain the proposed algorithm to compute a solution of the governing equation~\eqref{eq:EkEq_Lap} on a polyhedron mesh. In Section~\ref{sec:num}, numerical properties of the propose algorithm are presented by examples with exact solutions. Finally, we conclude in Section~\ref{sec:conclusion}

\section{Proposed method}\label{sec:prop_alg}

We start with explaining concrete notations to bring a better understanding of polyhedral cells. In the following subsections, a linearized eikonal equation with Laplacian regularizer is introduced and its discretization based on a cell-centered FVM is presented in details. Finally, we explain how to design a decreasing sequence of regularization parameters and propose an algorithm to approximate a viscosity solution of~\eqref{eq:EkEq_org} by solving~\eqref{eq:EkEq_Lap} in the last subsection. 

\subsection{Notations}
\begin{figure}
	\begin{center}
		\begin{tabular}{c}
			\includegraphics[height=2.5cm]{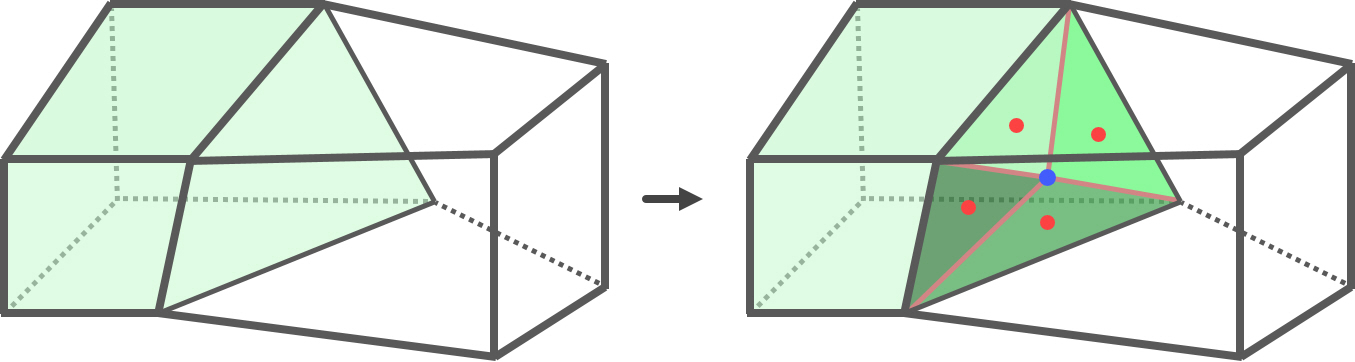}
		\end{tabular}
	\end{center}
	\caption{An illustration of two polyhedral cells with a tessellated face.}\label{fig:cell_dig}
\end{figure}

Let us denote a discretized computational domain as a union of non-overlapped polyhedral cells with a non-zero volume:
\begin{align}\label{eq:dis_dom}
	\bar{\Omega} = \bigcup_{p \in \mathcal{I}} \bar{\Omega}_p,
\end{align}
where $\Omega_p$ is open and $\mathcal{I}$ is a set of the indices of cells; see an illustration of two polyhedral cells in Figure~\ref{fig:cell_dig}. If a face is in-between two adjacent cells, we call it an internal face. Otherwise, we call it a boundary face. A set $\mathcal{G}$ is the collection of indices of all internal faces. Since a face of a polyhedron cell is difficult to be a plane in a general shape of computational domain, we always consider a tessellation of a face into triangles unless the face is already a triangle. From a face $e_g$, $g\in\mathcal{G}$, whose vertices are $\mathbf{x}_{v_i}$, $i=1,\ldots,r_g$, a triangle $\mathcal{T}_i = \mathcal{T}(\mathbf{x}_{v_i},\mathbf{x}_{v_{i+1}},\mathbf{x}_0)$ of three points, $\mathbf{x}_{v_i}$, $\mathbf{x}_{v_{i+1}}$, and the center of the mass $\mathbf{x}_{0} = \frac{1}{r_g}\sum_{i=1}^{r_g} \mathbf{x}_{v_i}$ is used to define a center of the face:
\begin{align}\label{eq:poly_face_ct}
	\mathbf{x}_g = \frac{\sum_{i=1}^{r_g} \left|\mathcal{T}_i\right| \bar{\mathbf{x}}_i }{\sum_{i=1}^{r_g} \left|\mathcal{T}_i\right|},
\end{align}
where $\mathbf{x}_{v_{r_g+1}} = \mathbf{x}_{v_1}$ and $\bar{\mathbf{x}}_i$ and $|\mathcal{T}_i|$ are the center and area of the triangle $\mathcal{T}_i$, respectively. Note that $\mathbf{x}_g$ is not necessarily same as $\mathbf{x}_{0}$ the center of the mass in general. In order to indicate the tessellated faces of a general face indexed by $\mathcal{G}$, we define a set of the indices of a tessellated internal and boundary faces as $\mathcal{F}$ and $\mathcal{B}$. For example, $e_f$, $f \in \mathcal{F}$, is a triangle on a face between left and right cells in Figure~\ref{fig:cell_dig} and $\mathbf{x}_f$ (red point) is the center of the triangle, where all triangles share a vertex, the center of the face $\mathbf{x}_g$ (blue point). To sum up, for a face $e_g$, $g\in\mathcal{G}$, there exists a subset $\mathcal{F}_g \subset \mathcal{F}$ such that 
\begin{align*}
	e_g = \bigcup_{f \in \mathcal{F}_g} e_f.
\end{align*}
If a face $e_g$ is not a triangle, it is a collection of tessellated faces (triangles) $e_f$, $f \in \mathcal{F}_g$, whose common vertex is $\mathbf{x}_g$. If $e_g$ is a triangle, then there is an index $f\in\mathcal{F}$ such that $e_g = e_f$.

For a cell $\Omega_p$, $p \in \mathcal{I}$, we define a set $\mathcal{N}_p$ as the indices of neighbor cells $\Omega_q$ such that the intersection $\partial \Omega_p \cap \partial \Omega_q = e_g$, $g\in\mathcal{G}$, is a face of non-zero area between two adjacent cells. We also define $\mathcal{F}_p$ and $\mathcal{B}_p$ as internal and boundary triangles tessellated by faces of $\Omega_p$. When $\mathcal{B}_p$ is empty, we call the cell $\Omega_p$ as an internal cell. Otherwise, it is called as a boundary cell. For example, if a green cell $\Omega_p$ in Figure~\ref{fig:cell_dig} is a boundary cell whose only left side is a part of the boundary of the computational domain, $|\mathcal{N}_p| = 5$, $|\mathcal{F}_p| = 20$, and $|\mathcal{B}_p| = 4$. If the cell next to the green cell is $\Omega_q$, then $q\in\mathcal{N}_p$ and there is an index $g \in \mathcal{G}$ such that $e_g = \partial \Omega_p \cap \partial \Omega_q$. In the rest of paper, we use the subscripts $f$, $b$, and $g$ to indicate an internal triangle $e_f$, a boundary triangle $e_b$, and an internal face $e_g$, respectively, unless otherwise noted.

For an internal triangle $e_f$, $f\in\mathcal{F}_p$, $p\in\mathcal{I}$, the vector $\mathbf{n}_{pf}$ is the outward normal to the triangle and its length is the area of the triangle, $|\mathbf{n}_{pf}| = |e_f|$. Then, $e_f \subset \partial \Omega_q$ for $q\in\mathcal{N}_p$, $\mathbf{n}_{qf} = -\mathbf{n}_{pf}$ holds. For a boundary triangle $e_b$, $b\in\mathcal{B}_p$, $p\in\mathcal{I}$, the vector $\mathbf{n}_{b} = \mathbf{n}_{pb}$ is the outward normal to the triangle, that is, the outward normal to the boundary of the computational domain, and its length is the area of the triangle, $|\mathbf{n}_{b}| = |e_b|$. When a directional vector is specified by two position vectors $\mathbf{x}_a$ and $\mathbf{x}_b$, we use a notation $\mathbf{d}_{ab} = \mathbf{x}_b - \mathbf{x}_a$. For an internal face $e_g=\partial \Omega_p \cap \partial \Omega_q$, $g\in\mathcal{G}$, $p,\:q\in\mathcal{I}$, whose vertices are written by $\mathbf{x}_{v_i}$, $i=1,\ldots,r_g$, we define a vector:
\begin{align}\label{eq:poly_face_nor}
	\mathbf{n}_g = \frac{1}{2}\sum_{i=2}^{r_g-1} \mathbf{d}_{v_1 v_i} \times \mathbf{d}_{v_1 v_{i+1}},
\end{align}
where the order of vertices is decided such that the cross product $\mathbf{d}_{v_1 v_i} \times \mathbf{d}_{v_1 v_{i+1}}$ indicates the outward to the cell $\Omega_p$ for all $i=2,\ldots,r_g-1$. If the face $e_g$ is planar, the vector $\mathbf{n}_g$ becomes an outward normal vector to the face of the cell $\Omega_p$ and its length $|\mathbf{n}_g| = |e_g|$ is the area of the face.

The characteristic length of a discretized domain $\cup_{p\in\mathcal{I}_{\text{L}}} \Omega_p$ is defined by the average of one-third power to the volume of the bounding box of a cell:
\begin{align}\label{eq:char_len}
	h_{\text{L}} = \frac{1}{|\mathcal{I}_{\text{L}}|} \sum_{p \in \mathcal{I}_{\text{L}}} |\Omega_p|_{B}^{\frac{1}{3}},
\end{align}
where $|\Omega_p|_{B}$ is the volume of the box whose diagonal is a vector $\mathbf{x}_M - \mathbf{x}_m$, where $\mathbf{x}_m$ and $\mathbf{x}_M$ are componentwise minimum and maximum of all points $\mathbf{x}_{v_i}$, for $i=1,\ldots,r_g$, respectively. The L indicates the level of mesh refinement, that is, when L increases, finer cells are generated. In Section~\ref{sec:num}, we use four levels of cells, roughly $h_{{\text{L}} + 1} \approx \frac{1}{2} h_{\text{L}}$, to check the experimental order of convergence ($EOC$).

\subsection{Linearized eikonal equation with Laplacian regularization}
In this subsection, we assume that there is a known function $u_{\epsilon'}$ which is possibly close to the solution of~\eqref{eq:EkEq_Lap} with a regularization parameter $\epsilon' > 0$. We present how to use a cell-centered finite volume method with the Soner boundary condition to numerically find a solution $u_{\epsilon}$ of~\eqref{eq:EkEq_Lap} with a smaller regularization parameter $\epsilon < \epsilon'$. Firstly, a linearization of the nonlinear term in~\eqref{eq:EkEq_Lap} is used to obtain an equation of unknown function $u_\epsilon$:
\begin{equation}\label{eq:EkEq_Lin}
	-\epsilon \triangle u_\epsilon(\mathbf{x}) + \mathbf{v}(\mathbf{x})\cdot \nabla u_\epsilon(\mathbf{x}) = 1,  \quad
	\mathbf{v}(\mathbf{x}) = \frac{\nabla u_{\epsilon'}(\mathbf{x})}{|\nabla u_{\epsilon'}(\mathbf{x})|_{\sigma}}, \quad \mathbf{x} \in \Omega\setminus\Gamma,
\end{equation}
where $|\mathbf{x}|_\sigma = (|\mathbf{x}|^2 + \sigma^2)^{\frac{1}{2}}$ with a small constant $\sigma = 10^{-12}$. Note that $\mathbf{v}$ is a fixed vector and the details of computing $u_{\epsilon'}$ is explained in the next subsection. Secondly, we show how to apply the Soner boundary condition in a cell-centered finite volume method. Even if a discretization of the normal flow term, $\mathbf{v} \cdot \nabla u_{\epsilon}$, with Soner boundary condition is already presented in~\cite{ref:HMFB21}, we repeat the key points of the numerical scheme in order to completely explain a discretization of the Laplacian term with Soner boundary condition based on the flux-balanced approximation~\cite{ref:FMHMB19} on a polyhedral cell.

Before we derive a discretization of using Soner boundary condition, a gradient computation is necessary at the center $\mathbf{x}_p$ of the cell $\Omega_p$. Since $\Gamma$ can be a part of the boundary of the computational domain, let us denote an index set to indicate triangles on the boundary and $\Gamma$:
\begin{align}\label{eq:bdry_D}
	\mathcal{B}_D = \{ b \in \mathcal{B} : e_b \subset \Gamma \cap \partial \Omega\}.
\end{align}
Defining $\mathcal{A}_p = \mathcal{N}_p \cup \left( \mathcal{B}_p \cap \mathcal{B}_D \right)$, the weighted least-squares method is used to compute the gradient at the center $\mathbf{x}_p$:
\begin{align}\label{eq:grad_p}
	\nabla u_p \equiv \nabla u(\mathbf{x}_p) = \argmin \limits_{ \substack{\mathbf{y} \in \mathbb{R}^3 \\ |\mathbf{y}| \leq 1}
	} \left( \sum_{a \in \mathcal{A}_p} \frac{(u_p + \mathbf{y} \cdot \mathbf{d}_{pa} - u_a)^2}{|\mathbf{d}_{pa}|^2} \right).
\end{align}
Note that $u_a = u(\mathbf{x}_a) = 0$, $a \in \mathcal{B}_p \cap \mathcal{B}_D$, because of Dirichlet boundary condition in~\eqref{eq:EkEq_Lap}. The constraint in~\eqref{eq:grad_p} is also used in~\cite{ref:HMFB21} which brings a more stable numerical computation. A componentwise constraint of the gradient is presented in~\cite{ref:TRSBB05,ref:XT10} to improve a stability. In~\cite{ref:FB20}, the same constraint in~\eqref{eq:grad_p} is shown for a variational approach to solve the eikonal equation.

Now, we use the basic idea of flux-balanced approximation~\cite{ref:FMHMB19} and a deferred correction method with a concept of inflow-implicit outflow-explicit method on the linearized equation~\eqref{eq:EkEq_Lin}. By the relation $\nabla u \cdot \mathbf{v}= \nabla \cdot (u \mathbf{v}) - u \nabla \cdot \mathbf{v}$, the equation is evaluated at the center of the cell $\Omega_p$:
\begin{align*}
	-\epsilon \nabla \cdot \nabla u(\mathbf{x}_p) + \nabla \cdot (u \mathbf{v} )(\mathbf{x}_p) - u(\mathbf{x}_p) \nabla \cdot \mathbf{v} (\mathbf{x}_p) = 1,
\end{align*}
where $u = u_\epsilon$ for simplicity of formula derivation. Approximating a divergence of vector-valued function $\mathbf{F}$ evaluated at $\mathbf{x}_p$ by integrating over the cell $\Omega_p$:
\begin{align*}
	\nabla \cdot \mathbf{F} (\mathbf{x}_p) \approx \frac{1}{|\Omega_p|} \int_{\Omega_p} \nabla \cdot \mathbf{F} dV = \frac{1}{|\Omega_p|} \int_{\partial \Omega_p} \mathbf{F} \cdot \mathbf{n} dS,
\end{align*}
where $\mathbf{n}$ is an unit outward normal vector to $\partial \Omega_p$, then we have 
\begin{align}\label{eq:dis_form}
	\begin{split}
		0 &= -\epsilon \int_{\partial \Omega_p} \nabla u \cdot \mathbf{n} dS + \int_{\partial \Omega_p} u \mathbf{v} \cdot \mathbf{n} dS - u_p \int_{\partial \Omega_p} \mathbf{v} \cdot \mathbf{n} dS - |\Omega_p| \\
		&= -\epsilon \texttt{(II)} + \left( \texttt{(I)} - |\Omega_p| \right) 
	\end{split}
\end{align} 
After the complete discretization of two terms $\texttt{(I)}$ and $\texttt{(II)}$ is derived, we present a deferred correction method to compute the solution of~\eqref{eq:EkEq_Lin} in the end of this subsection.

The term $\texttt{(I)}$ in~\eqref{eq:dis_form} is further calculated:
\begin{align}\label{eq:FirstTerm_int}
	\texttt{(I)} &= \sum_{f \in \mathcal{F}_p \cup \mathcal{B}_p} \left( \int_{e_f} u \mathbf{v} \cdot \frac{\mathbf{n}_{pf}}{|\mathbf{n}_{pf}|} dS- u_p \int_{e_f} \mathbf{v} \cdot \frac{\mathbf{n}_{pf}}{|\mathbf{n}_{pf}|} dS \right) \approx \sum_{f \in \mathcal{F}_p \cup \mathcal{B}_p} \left( u_{pf} - u_p \right) \mu_{pf},
\end{align}
where $u_{pf}$ is a value at the center of face $e_f$, $f\in \mathcal{F}_p$, $u_p = u(\mathbf{x}_p)$, and the normal flux $\mu_{pf}$ is computed by
\begin{align}
	\mu_{pf} = \int_{e_f} \mathbf{v} \cdot \frac{\mathbf{n}_{pf}}{|\mathbf{n}_{pf}|} dS \approx \mathbf{v}_{f} \cdot \mathbf{n}_{pf}.
\end{align}
The last term above is obtained by a formula with a small constant $\sigma = 10^{-12}$:
\begin{align}
	\mu_{pf} \approx \frac{\bm{\beta}_f}{\left(|\bm{\beta}_f|^2 + \sigma^2\right)^{\frac{1}{2}}} \cdot \mathbf{n}_{pf},
\end{align}
where $\bm{\beta}_f$ is a gradient whose length is less than $1$ at the center of the triangle $e_f$, $f \in \mathcal{F}_p \cup \mathcal{B}_p$, computed by a constraint minimization using pre-computed known function $u_{\epsilon'}$; see the equation (33) and the \textit{Remark 1} in~\cite{ref:HMFB21} for the technical details. In order to find the complete discretization of the first term, we define the sets of indices depending on the sign of the normal flux:
\begin{align}\label{eq:disjoints}
	\begin{split}
		\mathcal{B}^{-}_{p} = \{ b \in \mathcal{B}_{p} : \mu_{pb} < 0 \}, &\quad \mathcal{B}^{+}_{p} = \mathcal{B}_{p} \setminus \mathcal{B}^{-}_{p}, \\
		\mathcal{F}^{-}_{p} = \{ f \in \mathcal{F}_{p} : \mu_{pf} < 0 \}, &\quad \mathcal{F}^{+}_{p} = \mathcal{F}_{p} \setminus \mathcal{F}^{-}_{p}.
	\end{split}
\end{align}
Considering a general case of $\Gamma$ in~\eqref{eq:EkEq_Lap}, for example a part of $\partial \Omega$, we split the index set of boundary triangles into three cases:
\begin{align}\label{eq:split_bdry_nor}
	\mathcal{B}_p = \left( \mathcal{B}_p^{-} \cap \mathcal{B}_D \right) \cup \left( \mathcal{B}_p^{-} \setminus \mathcal{B}_D \right) \cup \mathcal{B}_p^{+},
\end{align}
where $\mathcal{B}_D$ is defined by~\eqref{eq:bdry_D}. On a boundary triangle $e_b$, $b \in \mathcal{B}_p$, we derive the numerical scheme on $\mathcal{B}_p^{+}$ because it does not violate Soner boundary condition and on $ \mathcal{B}_p^{-} \cap \mathcal{B}_D$ because Dirichlet boundary condition should be explicitly applied. The terms occurring on $\mathcal{B}_p^{-} \setminus \mathcal{B}_D$ should be set to zero not to violate the Soner boundary condition. 
Then, the original discretization of the normal flow in~\cite{ref:HMFB17_1,ref:HMFMB19} is changed because of using the Soner boundary condition:
\begin{align}\label{eq:FirstTerm_dis}
	\begin{split}
		\texttt{(I)} &\approx \sum_{f \in \mathcal{F}^{-}_{p}} \left( u_q + \mathcal{D}_q u \cdot \mathbf{d}_{qf}  - u_p \right) \mu_{pf} + \sum_{f \in \mathcal{B}^{+}_{p} \cup \mathcal{F}^{+}_{p} } \left( \mathcal{D}_p u \cdot \mathbf{d}_{pf} \right) \mu_{pf} \\
		&+ \sum_{b \in \mathcal{B}^{-}_{p} \cap \mathcal{B}_{D} } \left( u_b - u_p \right)\mu_{pb}
	\end{split}
\end{align}
where $u_b = u(\mathbf{x}_b) = 0$, $b \in \mathcal{B}^{-}_p \cap \mathcal{B}_D$ and the modified inflow-based gradient $\mathcal{D}_p u$ is used to include the influence of the Soner boundary condition:
\begin{align}\label{eq:mod_IBG}
	\mathcal{D}_p u = 
	\displaystyle \frac{\displaystyle \sum_{f \in \mathcal{F}_p^{-} \cup (\mathcal{B}^{-}_p \cap \mathcal{B}_D)} \frac{1}{|\mathbf{d}_{pf}|} \bm{\beta}_f }{\displaystyle \sum_{f \in \mathcal{F}_p^{-} \cup (\mathcal{B}^{-}_p \cap \mathcal{B}_D) } \frac{1}{|\mathbf{d}_{pf}|}}.
\end{align}

The term $\texttt{(II)}$ is followed by the discretization of flux-balanced approximation~\cite{ref:FMHMB19}:
\begin{align}\label{eq:SecTerm_int}
	\texttt{(II)} = \sum_{q \in \mathcal{N}_p} \int_{e_g} \nabla u \cdot \frac{\mathbf{n}_{g}}{|\mathbf{n}_{g}|} dS + \sum_{b \in \mathcal{B}_p} \int_{e_b} \nabla u \cdot \frac{\mathbf{n}_{pb}}{|\mathbf{n}_{pb}|} dS,
\end{align}
where a polygon face $e_g = \partial \Omega_p \cap \partial \Omega_q$, $g \in \mathcal{G}$, $q \in \mathcal{N}_p$, $p \in \mathcal{I}$. From the centers of two cells, $\mathbf{x}_p$ and $\mathbf{x}_q$, we find two points $\mathbf{x}_{p'}$ and $\mathbf{x}_{q'}$ such that the directional vectors $\mathbf{d}_{pp'}$ and $\mathbf{d}_{qq'}$ are perpendicular to the line passing at $\mathbf{x}_g$~\eqref{eq:poly_face_ct} along the direction $\mathbf{n}_g$~\eqref{eq:poly_face_nor}:
\begin{align*}
	\mathbf{d}_{pp'}  = \mathbf{d}_{pg} - \left(\frac{\mathbf{n}_{g}}{|\mathbf{n}_{g}|} \cdot \mathbf{d}_{pg}\right) \frac{\mathbf{n}_{g}}{|\mathbf{n}_{g}|}, \quad
	\mathbf{d}_{qq'}  = \mathbf{d}_{qg} - \left(\frac{\mathbf{n}_{g}}{|\mathbf{n}_{g}|} \cdot \mathbf{d}_{qg}\right) \frac{\mathbf{n}_{g}}{|\mathbf{n}_{g}|}.
\end{align*}
Using the explicit expression $\mathbf{d}_{pp'}$ and $\mathbf{d}_{qq'}$, we have an approximation of the first integral in~\eqref{eq:SecTerm_int}:
\begin{align}\label{eq:SecTerm_dis_1}
	\begin{split}
		\sum_{q \in \mathcal{N}_p} \int_{e_g} \nabla u \cdot \frac{\mathbf{n}_{g}}{|\mathbf{n}_{g}|} dS &\approx \sum_{q \in \mathcal{N}_p} \frac{|e_g|}{|\mathbf{d}_{p'q'}|} \left( u_{q'} -  u_{p'}\right) \\
		&\approx \sum_{q \in \mathcal{N}_p} \frac{|e_g|}{|\mathbf{d}_{p'q'}|} \left( \left(  u_q + \nabla u_q \cdot \mathbf{d}_{qq'} \right) - \left( u_p + \nabla u_p \cdot \mathbf{d}_{pp'} \right) \right)
	\end{split}
\end{align}
Note that more technical details are described in~\cite{ref:FMHMB19}. The second integral in~\eqref{eq:SecTerm_int} should be considered more carefully to apply the Soner boundary condition. Similar to~\eqref{eq:split_bdry_nor}, we split the index set of $\mathcal{B}_p$ into three cases.
\begin{align}\label{eq:split_bdry_Lap}
	\mathcal{B}_p = \left( \mathcal{B}_p \cap \mathcal{B}_D \right) \cup \left( \mathcal{B}_p \setminus \mathcal{B}_D \right) = \left( \mathcal{B}_p \cap \mathcal{B}_D \right) \cup \left( \mathcal{B}^{+}_p \setminus \mathcal{B}_D \right) \cup \left( \mathcal{B}^{-}_p \setminus \mathcal{B}_D \right),
\end{align}
The first case, on a triangle $e_b$, $b\in\mathcal{B}_p \cap \mathcal{B}_D$, the Dirichlet condition is applied. The second case, we use numerical values inside the computational domain. The third case, the terms occurring on $\mathcal{B}_p^{-} \setminus \mathcal{B}_D$ should be set to zero not to violate the Soner boundary condition. Considering the mentioned three cases, we have an approximation of the second integral~\eqref{eq:SecTerm_int}:
\begin{align}\label{eq:SecTerm_dis_2}
	\sum_{b \in \mathcal{B}_p} \int_{e_b} \nabla u \cdot \frac{\mathbf{n}_{pb}}{|\mathbf{n}_{pb}|} dS \approx \sum_{b \in \mathcal{B}_p \cap \mathcal{B}_D}  \frac{|e_b|}{|\mathbf{d}_{p'b}|} \left( u_b -  u_p - \nabla u_p \cdot \mathbf{d}_{pp'} \right) + \sum_{b \in \mathcal{B}_p^{+} \setminus \mathcal{B}_D} \nabla  u_p \cdot \mathbf{n}_{b},
\end{align}
where $u_b = u(\mathbf{x}_b) = 0$, $b \in \mathcal{B}_p \cap \mathcal{B}_D$.

Combining all derivations~\eqref{eq:FirstTerm_dis},~\eqref{eq:SecTerm_dis_1}, and~\eqref{eq:SecTerm_dis_2}, we have a complete discretization using the Soner boundary condition to solve~\eqref{eq:EkEq_Lin}:
\begin{align}\label{eq:EkEq_Lin_sys}
	\begin{split}
		0 = & - \epsilon \left( \sum_{q \in \mathcal{N}_p} \frac{|e_g|}{|\mathbf{d}_{p'q'}|} \left(  u_q + \nabla u_q \cdot \mathbf{d}_{qq'} - u_p - \nabla u_p \cdot \mathbf{d}_{pp'} \right) \right) \\
		& - \epsilon \left( \sum_{b \in \mathcal{B}_p \cap \mathcal{B}_D}  \frac{|e_b|}{|\mathbf{d}_{p'b}|} \left( u_b -  u_p - \nabla u_p \cdot \mathbf{d}_{pp'} \right) + \sum_{b \in \mathcal{B}_p^{+} \setminus \mathcal{B}_D} \nabla  u_p \cdot \mathbf{n}_{b} \right) \\
		&+ \sum_{f \in \mathcal{F}^{-}_{p}} \left( u_q + \mathcal{D}_q u \cdot \mathbf{d}_{qf}  - u_p \right) \mu_{pf} + \sum_{f \in \mathcal{B}^{+}_{p} \cup \mathcal{F}^{+}_{p} } \left( \mathcal{D}_p u \cdot \mathbf{d}_{pf} \right) \mu_{pf} \\
		&+ \sum_{b \in \mathcal{B}^{-}_{p} \cap \mathcal{B}_{D} } \left( u_b - u_p \right)\mu_{pb} - |\Omega_p|,
	\end{split}
\end{align}
where the gradient $\nabla u_p$ and the modified inflow-based gradient $\mathcal{D}_p u$ are defined by~\eqref{eq:grad_p} and~\eqref{eq:mod_IBG}, respectively. On a regular cubic mesh, the displacement $\mathbf{d}_{pp'}$ and $\mathbf{d}_{qq'}$ are zero vectors and the equation above is a band block diagonal matrix equation. In parallel computing using domain decomposition with the $1$-ring neighborhood structure, if $\Omega_p$ is located in the  domain $D_1$ and one of faces of $\Omega_p$ is located between two domains, $D_1$ and $D_2$, one of the second face neighbor cells on $\Omega_p$, that is, $\Omega_r$, $r\in\mathcal{N}_q\setminus\mathcal{N}_p$ and $q\in\mathcal{N}_p$, may not be accessible by the domain $D_1$ where $\Omega_p$ is located. Such a cell exists when we compute $\nabla u_q$ or $\mathcal{D}_q u$ in the formulation of~\eqref{eq:EkEq_Lin_sys} and then it is not possible to construct a correct linear system in the domains $D_1$ and $D_2$. To overcome the mentioned technical difficult, we use a deferred correction method~\cite{ref:BHS84} to solve~\eqref{eq:EkEq_Lin_sys} iteratively:
\begin{align}\label{eq:EkEq_Lin_iter}
	\begin{split}
		0 = & - \epsilon \left( \sum_{q \in \mathcal{N}_p} \frac{|e_g|}{|\mathbf{d}_{p'q'}|} \left(  u_q^{k} + \nabla u_q^{k-1} \cdot \mathbf{d}_{qq'} - u_p^{k} - \nabla u_p^{k-1} \cdot \mathbf{d}_{pp'} \right) \right) \\
		& - \epsilon \left( \sum_{b \in \mathcal{B}_p \cap \mathcal{B}_D}  \frac{|e_b|}{|\mathbf{d}_{p'b}|} \left( u_b -  u_p^{k} - \nabla u_p^{k-1} \cdot \mathbf{d}_{pp'} \right) + \sum_{b \in \mathcal{B}_p^{+} \setminus \mathcal{B}_D} \nabla  u_p^{k-1} \cdot \mathbf{n}_{b} \right) \\
		&+ \sum_{f \in \mathcal{F}^{-}_{p}} \left( u_q^{k} + \mathcal{D}_q^{k-1} u \cdot \mathbf{d}_{qf}  - u_p^{k} \right) \mu_{pf} + \sum_{f \in \mathcal{B}^{+}_{p} \cup \mathcal{F}^{+}_{p} } \left( \mathcal{D}_p u^{k-1} \cdot \mathbf{d}_{pf} \right) \mu_{pf} \\
		&+ \sum_{b \in \mathcal{B}^{-}_{p} \cap \mathcal{B}_{D} } \left( u_b - u_p^{k} \right)\mu_{pb} - |\Omega_p|,
	\end{split}
\end{align}
where $k\in\mathbb{N}$ and $u^{0} = u_{\epsilon'}$. Keeping in mind the formulation above, we continue to discuss a decreasing sequence of regularization parameters $\epsilon$ and a pre-computed function $u_{\epsilon'}$ in~\eqref{eq:EkEq_Lin} in the next subsection.

\subsection{The regularization parameter $\epsilon$}\label{sec:num:reg}

In this subsection, the proposed algorithm is described. Firstly, we explain two observations of the regularization parameter $\epsilon$ in numerical points of view. Secondly, considering the observations, we propose a sequential algorithm to solve~\eqref{eq:EkEq_Lap}.

The vanishing viscosity method expects that the solution $u_\epsilon$ of~\eqref{eq:EkEq_Lap} becomes close to the viscosity solution of~\eqref{eq:EkEq_org} when the regularization parameter $\epsilon > 0$ is smaller and smaller. Similarly, we would like to find a numerical solution of~\eqref{eq:EkEq_Lap} converges to the viscosity solution when the characteristic length $h_L$~\eqref{eq:char_len} becomes smaller and smaller. That is, a numerical convergence is related to not only the characteristic length~$h_L$ but also the regularization parameter~$\epsilon$. An empirical relation between $h_L$ and $\epsilon$ to obtain a numerical convergence is that when~$h_L$ becomes smaller, $\epsilon$ must become smaller too. Such a relation is also observed in solving a variant of the phase field model of the simplified Stefan problem~\cite{ref:SB13}.

Another aspect of regularization parameter $\epsilon$ is that it cannot be too small in a fixed discretized domain. The reason is similar to that the time step cannot be too large in the time-relaxed eikonal equation~\eqref{eq:EkEq_Time}. The direct effect of time relaxation in a linear system is to add positive values on a diagonal element which brings more stable computation to solve the linear system; see more details in~\cite{ref:HMFB21}. When the time step is too large, the positive value being added to the diagonal elements is too small and then we can observe oscillation over the time as it is already shown in~\cite{ref:HMFB21}. Similarly, if the regularization parameter~$\epsilon$ is too small on a fixed discretized domain, then the numerical solution does not become close enough to the viscosity solution of~\eqref{eq:EkEq_org}. The same phenomenon of a regularization parameter $\eta$ is also observed in~\cite{ref:CWW13,ref:CV19} by solving a singularly perturbed boundary value problem in~\cite{ref:V67},
\begin{equation}
	\begin{alignedat}{2}
		-\eta^2 \triangle w(\mathbf{x}) + w(\mathbf{x}) &= 0, \quad &&\mathbf{x} \in \Omega, \\
		w(\mathbf{x}) &= 1, \quad &&\mathbf{x} \in \partial \Omega,
	\end{alignedat}
\end{equation}
which can be transformed to
\begin{equation}
	\begin{alignedat}{2}
		-\eta \triangle v(\mathbf{x}) + |\nabla v(\mathbf{x})|^2 &= 1, \quad &&\mathbf{x} \in \Omega, \\
		v(\mathbf{x}) &= 0, \quad &&\mathbf{x} \in \partial \Omega,
	\end{alignedat}
\end{equation}
by the Hopf–Cole transformation~\cite{ref:E98,ref:GR09}. 

The obvious effect of using regularization parameter $\epsilon$ is to eliminate singularities and compute a smooth solution. However, if the parameter is too large, the numerical solution is not accurate enough to be the distance function. If it is too small, the numerical computation is not stable enough. Therefore, a reasonable choice of the regularization parameter $\epsilon$ is from a large value to a small value in a certain range. Considering mentioned observations, we propose an algorithm to compute a sequential numerical solution on a polyhedron mesh with the characteristic length~$h_L$:
\begin{algorithm}[t] 
	\small
	\caption{A procedure to compute a numerical solution of~\eqref{eq:EkEq_Lap}} \label{alg:prop_alg}
	\begin{algorithmic}
		\Procedure{}{} 
		\State Initialization $u^0 = 0$.
		\State Set $n=1$ and $K_1 = 1$ .
		\State Find a solution $u^1 = u^{1,1}$ of~\eqref{eq:prop_alg} with $u^{1,0} = u^0 = 0$.
		\For{$n\gets2$ to $5$}
		\State Set $u^{n,0} = u^{n-1}$.
		\State Set $k=1$.
		\While{$\rho^{n,k} \geq \eta$}\Comment{See~\eqref{eq:residual}.}
		\State Find a solution $u^{n,k}$ of~\eqref{eq:prop_alg} with $u^{n,k-1}$.
		\State  $k \gets k+1$
		\EndWhile 
		\EndFor
		\EndProcedure		
	\end{algorithmic}
\end{algorithm}
\begin{align}\label{eq:prop_alg}
	\begin{split}
		0 = & - \epsilon_n \left( \sum_{q \in \mathcal{N}_p} \frac{|e_g|}{|\mathbf{d}_{p'q'}|} \left(  u_q^{n,k} + \nabla u_q^{n,k-1} \cdot \mathbf{d}_{qq'} - u_p^{n,k} - \nabla u_p^{n,k-1} \cdot \mathbf{d}_{pp'} \right) \right) \\
		& - \epsilon_n \left( \sum_{b \in \mathcal{B}_p \cap \mathcal{B}_D}  \frac{|e_b|}{|\mathbf{d}_{p'b}|} \left(u_b -  u_p^{n,k} - \nabla u_p^{n,k-1} \cdot \mathbf{d}_{pp'} \right) + \sum_{b \in \mathcal{B}_p^{+} \setminus \mathcal{B}_D} \nabla  u_p^{n,k-1} \cdot \mathbf{n}_{b} \right) \\
		&+ \sum_{f \in \mathcal{F}^{-}_{p}} \left( u_q^{n,k} + \mathcal{D}_q u^{n,k-1} \cdot \mathbf{d}_{qf}  - u_p^{n,k} \right) \mu_{pf}^{n-1} + \sum_{f \in \mathcal{B}^{+}_{p} \cup \mathcal{F}^{+}_{p} } \left( \mathcal{D}_p u^{n,k-1} \cdot \mathbf{d}_{pf} \right) \mu_{pf}^{n-1} \\		
		&+ \sum_{b \in \mathcal{B}^{-}_{p} \cap \mathcal{B}_{D} } \left( u_b - u_p^{n,k} \right)\mu_{pb}^{n-1} - |\Omega_p|, \quad k=1,\ldots,K_n
	\end{split}
\end{align}
where $u^0 = 0$, $u^{n,0} = u^{n-1}$ is a pre-computed solution of~\eqref{eq:EkEq_Lin}, $K_n$ is defined by~\eqref{eq:residual}, and we choose regularization parameters as a decreasing sequence:
\begin{align}
	\epsilon_n = \left( h_L \right)^{\frac{1}{2}n}, \quad n = 1,2,\ldots,5.
\end{align}
The solution $u^{n,k}$ is computed by $u^{n-1}$, the parameter $\epsilon_n$, and $k$ number of iterations in~\eqref{eq:prop_alg}. Note that we explain how to numerically implement Dirichlet boundary condition~\eqref{eq:EkEq_Lap} in the linear system~\eqref{eq:prop_alg} iat the end of this subsection. Rewriting \eqref{eq:prop_alg} formally as a matrix equation,
\begin{align}\label{eq:mateq}
	\mathbf{A}^{n-1} u^{n,k} = \mathbf{f}(u^{n,k-1}). 
\end{align}
an algebraic multigrid method is used to solve the above equation. The $k^{\text{th}}$ iteration is stopped at the smallest $K_n$ such that a residual error is smaller than a chosen error bound $\eta= 10^{-8}$:
\begin{align}\label{eq:residual}
	K_n = \min\left\{ k \in \mathbb{N} : \rho^{n,k} =\frac{1}{|\mathcal{I}|} \sum_{p \in \mathcal{I}} \left| \left(\mathbf{A}^{n-1}\phi^{n,k} - \mathbf{f}(\phi^{n,k})\right)_p \right| < \eta\right\}, \quad n \geq 2.
\end{align}
where the parenthesis above with a subscript $\left( \mathbf{r} \right)_p$ denotes the $p^{\text{th}}$ component of the vector $\mathbf{r}$. Then, we define $u^{n} \equiv u^{n,K_n}$ for $n \geq 2$. When $n=1$, we use $K_n = 1$. The proposed algorithm is also presented step by step in Algorithm~\ref{alg:prop_alg}. 
\begin{remark}\label{rmk:parallel}
	In the matrix of the linear system \eqref{eq:prop_alg} on the $p^{\text{th}}$ row, the diagonal element is the coefficient of $u_p^{n,k}$ and all off diagonal elements are the coefficients of $u_q^{n,k}$, $q\in \mathcal{N}_p$. It means the system only uses neighbor cells across faces of $\Omega_p$. Then, an implementation of~\eqref{eq:prop_alg} in a standard cell-centered FVM code is straightforwardly done for parallel computing using domain decomposition with $1$-ring face neighborhood.
\end{remark}
When $n=1$ in the proposed algorithm~\eqref{eq:prop_alg}, the linear system computes a solution of the equation below because all gradients are zero with the initial choice $u^{1,0} = u^0 = 0$:
\begin{equation}\label{eq:EkEq_Lap_n1}
	\begin{alignedat}{2}
		-\epsilon \triangle \bar{u}(\mathbf{x}) &= 1 \quad &&\mathbf{x} \in \Omega\setminus\Gamma, \\
		\bar{u}(\mathbf{x}) &= 0 \quad && \mathbf{x} \in \Gamma, \\
		\bm{\nu}(\mathbf{x}) \cdot \nabla \bar{u}(\mathbf{x}) & = 0 \quad && \mathbf{x} \in \partial \Omega \setminus\Gamma,
	\end{alignedat}
\end{equation}
The direction of $\nabla \bar{u}$ is same as the gradient of the viscosity solution in~\eqref{eq:EkEq_org} because their zero level set $\Gamma$ is identical. Then, for $n \geq 2$, the normalized vector~$\mathbf{v}$ in~\eqref{eq:EkEq_Lin}  on~$\Gamma$ is already same as the vector $\mathbf{v}$ computed by the viscosity solution of~\eqref{eq:EkEq_org}. In the case of $\Gamma = \partial \Omega$, the solution of Poisson equation~\eqref{eq:EkEq_Lap_n1} is also used to approximate a distance function on a close neighborhood of $\Gamma$ by a normalization scheme~\cite{ref:T98}. In~\cite{ref:AA14}, it is argued that there is a proximity in $L^2$ sense between the solution of~\eqref{eq:EkEq_Lap_n1} and the distance function from $\Gamma$.

In order to complete the description of the proposed algorithm, we need to explain how the boundary value on $\Gamma$ is implemented in a polyhedron mesh because $\Gamma \subset \bar{\Omega}$ is generally located on a given mesh. To do so, we define index sets to select the cells where $\Gamma$ is located in $\bar{\Omega}$:
\begin{align}
	\begin{split}
		\mathcal{I}^{1} &= \left\{ p \in \mathcal{I} : \bar{\Omega}_p \cap \Gamma \neq \emptyset, \: 
		\Gamma \subseteq \partial \Omega \right\}, \\
		\mathcal{I}^{2} &= \left\{ p \in \mathcal{I} : \bar{\Omega}_p \cap \Gamma \neq \emptyset, \:\Gamma \subsetneq  \Omega, \:\Gamma \cap \partial \Omega = \emptyset \right\}.
	\end{split}
\end{align}
If $\Gamma$ is a general shape, an octree search and point-in-cell algorithms are used to define the index sets above. Let us define a function $\mathcal{F} : K \subset \mathcal{I} \rightarrow \mathcal{I}$ by $\mathcal{F}(K) = \left\{ q \in \mathcal{I} : q \in \mathcal{N}_p, \: \forall p \in K \right\} \cup K$. Now, we use the set $\Gamma^0 = \mathcal{F}(\mathcal{F}(\mathcal{I}^{2})) \cup \mathcal{F}(\mathcal{I}^{1})$ and it is straightforward to compute the exact distance value from $\Gamma$ at all points $\mathbf{x}_p$, $p \in \Gamma^0$. An octree search algorithm can find a short list of potential elements in $\Gamma$ to compute the shortest distance from $\mathbf{x}_p$ to $\Gamma$ and it is efficient enough because all points $\mathbf{x}_p$, $p\in\Gamma^0$, are close to $\Gamma$. Then, the computed distance value on $\Gamma^0$ is used in the proposed algorithm. That is, on the $p^{\text{th}}$ row of the matrix~\eqref{eq:prop_alg}, $p\in\Gamma^0$, we use the value and make all relevant off-diagonal element of $\Omega_p$ to be zero in the matrix.

\section{Numerical results}\label{sec:num}

\begin{figure}
	\begin{center}
		\begin{tabular}{cc}
			\includegraphics[height=3.5cm]{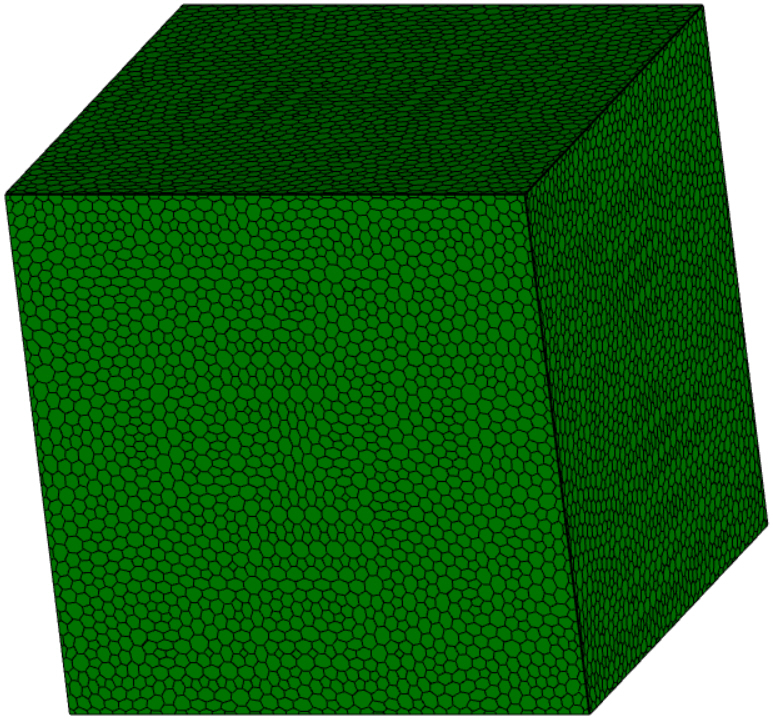} &
			\includegraphics[height=3.5cm]{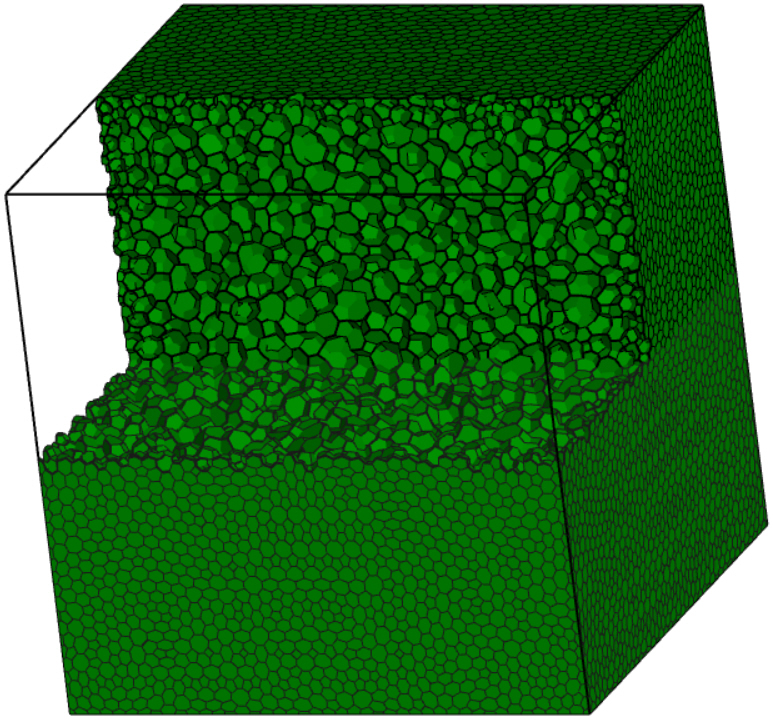} \\
			\multicolumn{2}{c}{(a) $\mathcal{M}^1_\text{L}$} \\ 
			\includegraphics[height=3.5cm]{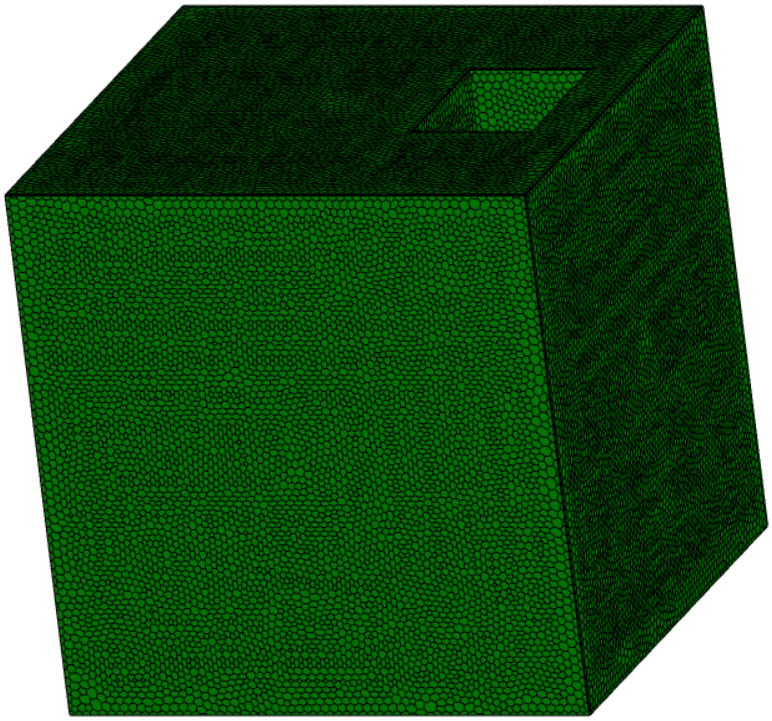} &
			\includegraphics[height=3.5cm]{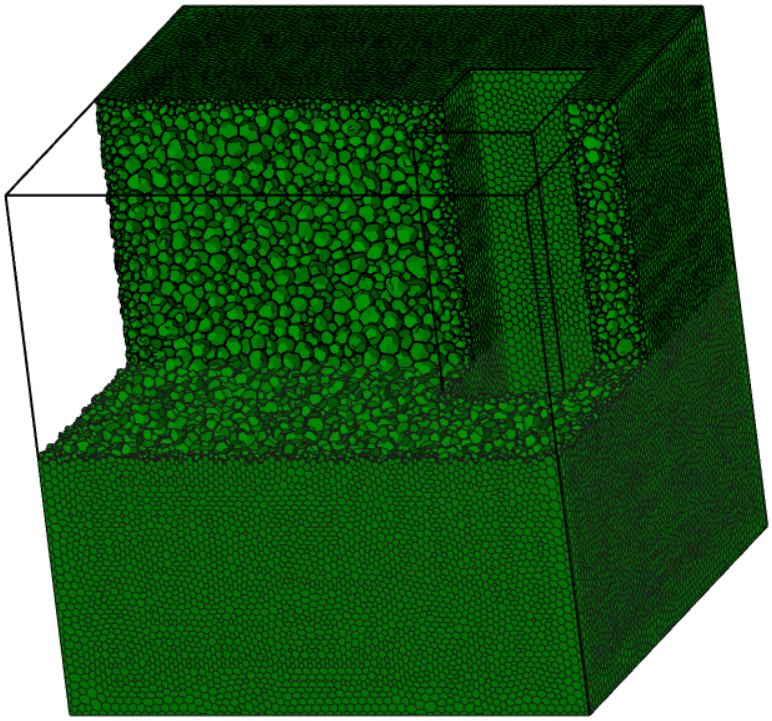} \\
			\multicolumn{2}{c}{(b) $\mathcal{M}^2_\text{L}$} \\ 
			\includegraphics[height=3.5cm]{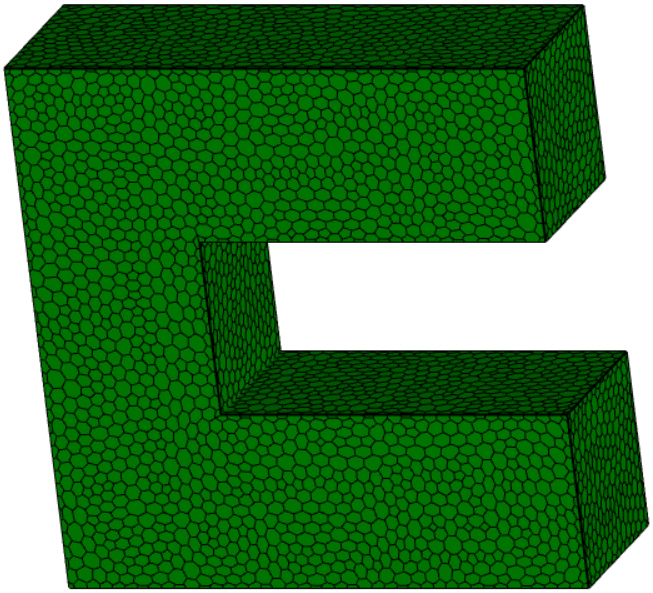} & 
			\includegraphics[height=3.5cm]{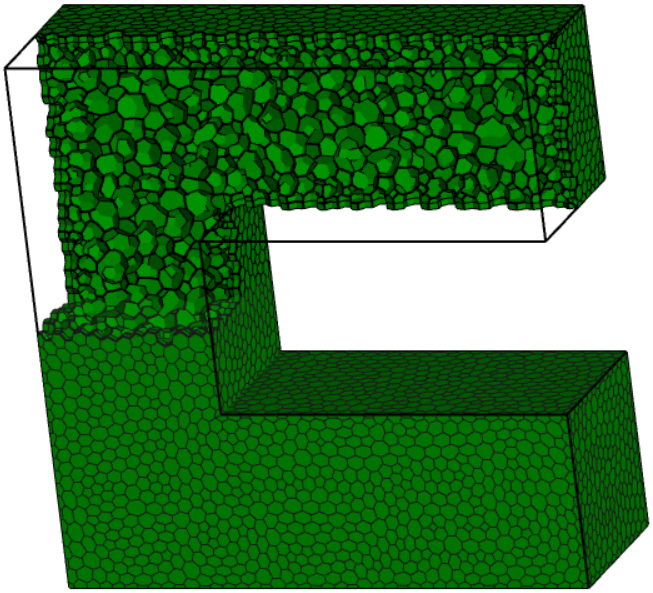} \\
			\multicolumn{2}{c}{(c) $\mathcal{M}^3_\text{L}$} \\ 
			\includegraphics[height=3.53cm]{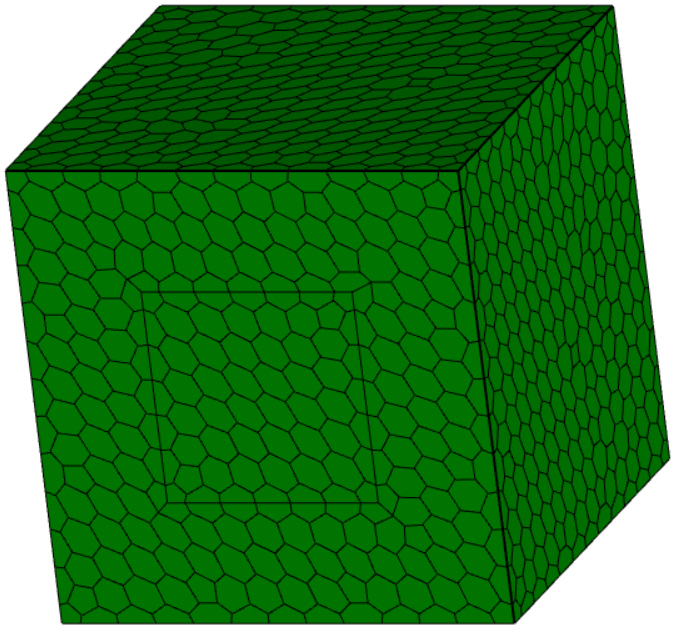} &
			\includegraphics[height=3.5cm]{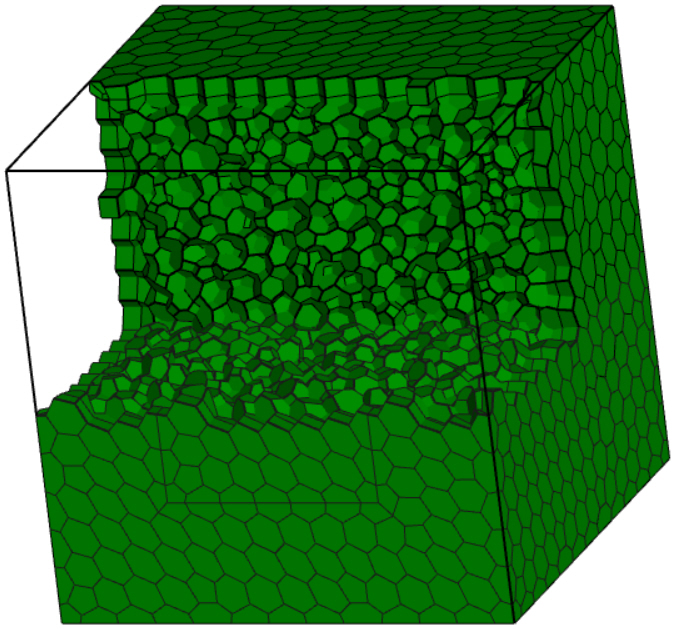} \\
			\multicolumn{2}{c}{(d) $\mathcal{M}^4_\text{L}$}
		\end{tabular}
	\end{center}
	\caption{It is an illustration of meshes for computational domains in Table~\ref{tab:meshes} with $L=1$. The left column shows the boundary of the computational domain and the right column shows what the polyhedral cells look like inside the domain. Note that the right side is the positive direction of $x$ axis, the top side is the positive direction of $y$ axis, and the direction coming out of the paper is the positive direction of $z$ axis.} \label{fig:meshes}
\end{figure}

\begin{table}[]
	\centering
	\small
	\begin{tabular}{
			c
			S[table-format=1]
			S[table-format=8]
			S[table-format=1.2,
			table-figures-exponent=2,
			table-sign-mantissa,
			table-sign-exponent]}
		\toprule
		mesh & \multicolumn{1}{c}{L} & \multicolumn{1}{c}{$|\mathcal{I}_{\text{L}}|$ } & \multicolumn{1}{c}{$h_{\text{L}}$} \\
		\midrule
		\multirow{4}{*}{$\mathcal{M}^1_{\text{L}}$} & 1 & 29954 & 9.91E-02 \\
		& 2 & 174917 & 5.63E-02 \\
		& 3 & 1151396 & 3.08E-02 \\
		& 4 & 8216986 & 1.63E-02 \\
		\cmidrule(lr){1-2} \cmidrule(lr){3-4}		
		\multirow{4}{*}{$\mathcal{M}^2_{\text{L}}$} & 1 & 129955 & 5.53E-02 \\   
		& 2 & 708104 & 3.22E-02 \\     
		& 3 & 4248440 & 1.84E-02 \\      
		& 4 & 28196165 & 1.02E-02 \\      
		\cmidrule(lr){1-2} \cmidrule(lr){3-4}
		\multirow{4}{*}{$\mathcal{M}^3_{\text{L}}$} & 1 & 18118 & 7.02E-02 \\
		& 2 & 74301 & 4.22E-02 \\
		& 3 & 362679 & 2.44E-02 \\
		& 4 & 1868820 & 1.45E-02 \\
		\cmidrule(lr){1-2} \cmidrule(lr){3-4}
		\multirow{4}{*}{$\mathcal{M}^4_{\text{L}}$} & 1 & 7863 & 6.52E-01 \\
		& 2 & 58091 & 3.42E-01 \\
		& 3 & 457436 & 1.73E-01 \\
		& 4 & 3660530 & 8.71E-02 \\		
		\bottomrule
	\end{tabular}
	\caption{The numbers of polyhedral cells $\mathcal{I}_{\text{L}}$ and the characteristic length $h_{\text{L}}$~\eqref{eq:char_len} of the meshes are presented; see the shape of the computational domains at the level $\text{L}=1$ in Figure~\ref{fig:meshes}.}\label{tab:meshes}
\end{table}

We present various examples to show numerical properties of the proposed algorithm~\eqref{eq:prop_alg}. The meshes generated by \texttt{AVL FIRE}\textsuperscript{\texttt{TM}} are illustrated in Figure~\ref{fig:meshes} and the number of polyhedral cells $|\mathcal{I}_{\text{L}}|$ and the characteristic length $h_{\text{L}}$~\eqref{eq:char_len} of the meshes are presented for four levels of meshes, $\text{L} \in \{1,2,3,4\}$, in Table~\ref{tab:meshes}. Note that $h_{\text{L+1}} < h_{\text{L}}$. The test examples are basically to compute a distance function from $\Gamma$ on a given computation domain $\Omega$ and all details are explained below with constants $\gamma_i = \frac{R_i}{15}$ for $i=1,\:2$, where $R_1=1.25$ and $R_2=10$.
\begin{enumerate}[leftmargin=3\parindent]
	\renewcommand{\labelenumi}{\textbf{\theenumi}}
	\renewcommand{\theenumi}{EX\arabic{enumi}}
	\item \label{ex_1} $\Gamma$ is a sphere with the center at the origin and the radius is $0.6$ in the computational domain $\Omega = [-R_1,R_1]^3$. The mesh $\mathcal{M}^1_{L}$ is used in Table~\ref{tab:meshes}. The first level of mesh is shown in Figure~\ref{fig:meshes}-(a).
	\item \label{ex_2} $\Gamma$ is a sphere at the origin with the radius is $0.3$ in the computational domain $\Omega = [-8\gamma_1,22\gamma_1]\times[-15\gamma_1,15\gamma_1]\times[-15\gamma_1,15\gamma_1] \setminus \Omega' $, where $\Omega' = [8\gamma_1, 15\gamma_1]\times[-15\gamma_1,15\gamma_1]\times[-5\gamma_1,5\gamma_1]$. The mesh $\mathcal{M}^2_{L}$ is used in Table~\ref{tab:meshes}. The first level of mesh is shown in Figure~\ref{fig:meshes}-(b). 
	\item \label{ex_3} The computational domain is $\Omega = [-15\gamma_1,15\gamma_1]\times[-15\gamma_1,15\gamma_1]\times[-5\gamma_1,5\gamma_1] \setminus \Omega'$, where $\Omega' = [-5\gamma_1,15\gamma_1]\times[-5\gamma_1,5\gamma_1]\times[-5\gamma_1,5\gamma_1]$ and $\Gamma = \{15\gamma_1\} \times[5\gamma_1,15\gamma_1] \times[-5\gamma_1,5\gamma_1]$ is the upper right plane. The mesh $\mathcal{M}^3_{L}$ is used in Table~\ref{tab:meshes}. The first level of mesh is shown in Figure~\ref{fig:meshes}-(c). 
	\item \label{ex_4} The computational domain $\Omega$ is same as \ref{ex_3} and $\Gamma = \{15\gamma_1\} \times[5\gamma_1,15\gamma_1] \times[-5\gamma_1,5\gamma_1] \cup \{15\gamma_1\} \times[-15\gamma_1,-5\gamma_1] \times[-5\gamma_1,5\gamma_1]$ is the upper right and the lower right planes in Figure~\ref{fig:meshes}-(c). The mesh $\mathcal{M}^3_{L}$ is used in Table~\ref{tab:meshes}.
	\item \label{ex_5} The computational domain $\Omega$ is same as \ref{ex_3} and $\Gamma =\partial \Omega$. The mesh $\mathcal{M}^3_{L}$ is used in Table~\ref{tab:meshes}.
	\item \label{ex_6} The computational domain is $\Omega = \left[-\frac{R_2}{2},\frac{R_2}{2}\right]^3$ and $\Gamma=\partial \Omega$. The mesh $\mathcal{M}^4_{L}$ is used in Table~\ref{tab:meshes}.
	\item \label{ex_7} The computational domain is same as  \ref{ex_6} and $\Gamma$ is a circle with the center at the origin and the radius $0.6$, where the normal vector of the plane containing the circle is $z$ axis. The mesh $\mathcal{M}^4_{L}$ is used in Table~\ref{tab:meshes}.
	\item \label{ex_8} The computational domain is same as  \ref{ex_6} and $\Gamma$ is a disk whose boundary is the circle in \ref{ex_7}. The mesh $\mathcal{M}^4_{L}$ is used in Table~\ref{tab:meshes}.
	\item \label{ex_9} The computational domain is $\Omega =  [-R_2,R_2]^3$ and $\Gamma$ is a square with the center at the origin and the length of the side is $7 r_2$, where the normal vector of the plane containing the square is $z$ axis. The mesh $\mathcal{M}^4_{L}$ is used in Table~\ref{tab:meshes}. The first level of mesh is shown in Figure~\ref{fig:meshes}-(d).
	\item \label{ex_10} The computational domain is same as  \ref{ex_9} and $\Gamma$ is two squares of the same size used in \ref{ex_9}. The center of the first and second square is located at $(0,0,7.5\gamma_2)$ and $(0,0,-7.5\gamma_2)$, respectively. The mesh $\mathcal{M}^4_{L}$ is used in Table~\ref{tab:meshes}.
\end{enumerate}
The examples from \ref{ex_1} to \ref{ex_5} are already presented in~\cite{ref:HMFB21} where one can find exact solutions of the examples explicitly and the exact solutions from \ref{ex_6} to \ref{ex_10} can be easily obtained. Note that the polyhedron mesh $\mathcal{M}^3$ is exactly same as the one in~\cite{ref:HMFB21}, but we use polyhedral meshes $\mathcal{M}^1$ and $\mathcal{M}^2$ whose characteristic length is slightly less than twice as small in~\cite{ref:HMFB21}. A typical body-fitted surface mesh is used on two squares $\Gamma$ in the case of~\ref{ex_10}. In Figure~\ref{fig:meshes}-(d), one of squares is visible on the boundary of $\mathcal{M}_L^4$, $L=1$. The same mesh is used to test cases from~\ref{ex_6} to~\ref{ex_10}.

\begin{figure}
	\begin{center}
		\begin{tabular}{cc}
			\includegraphics[height=3cm]{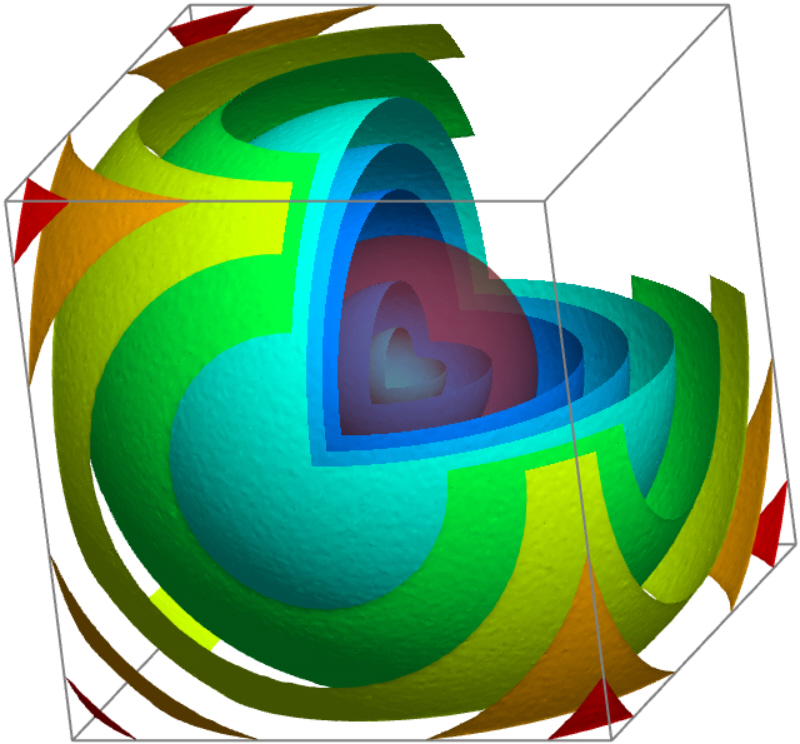} &
			\includegraphics[height=3cm]{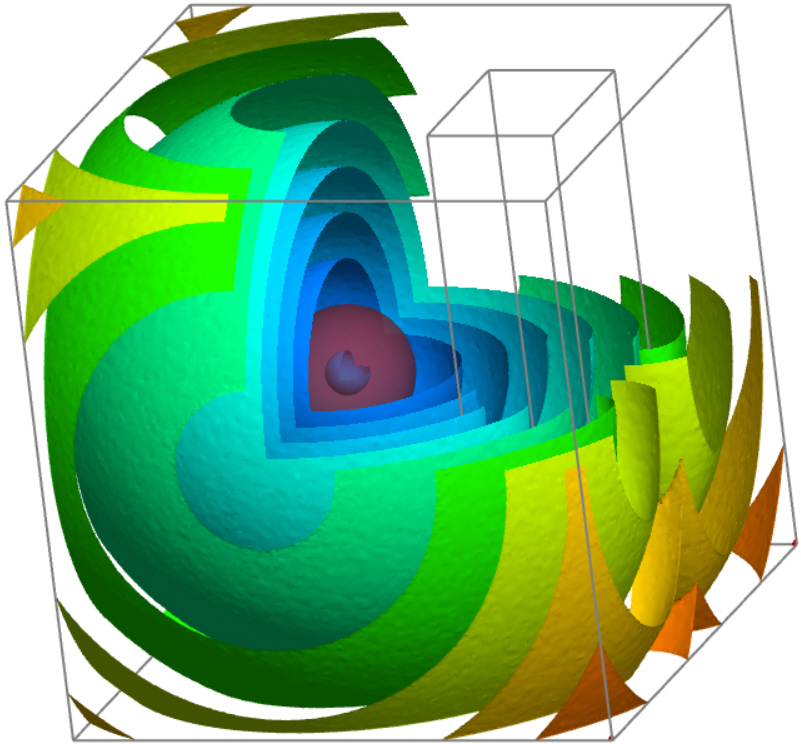} \\
			\ref{ex_1} & \ref{ex_2} \\
			\includegraphics[height=3cm]{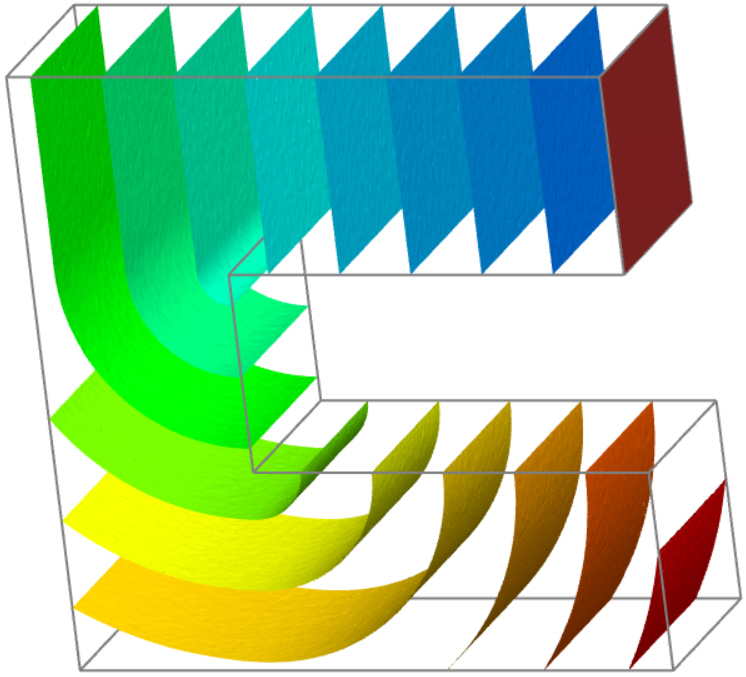} &
			\includegraphics[height=3cm]{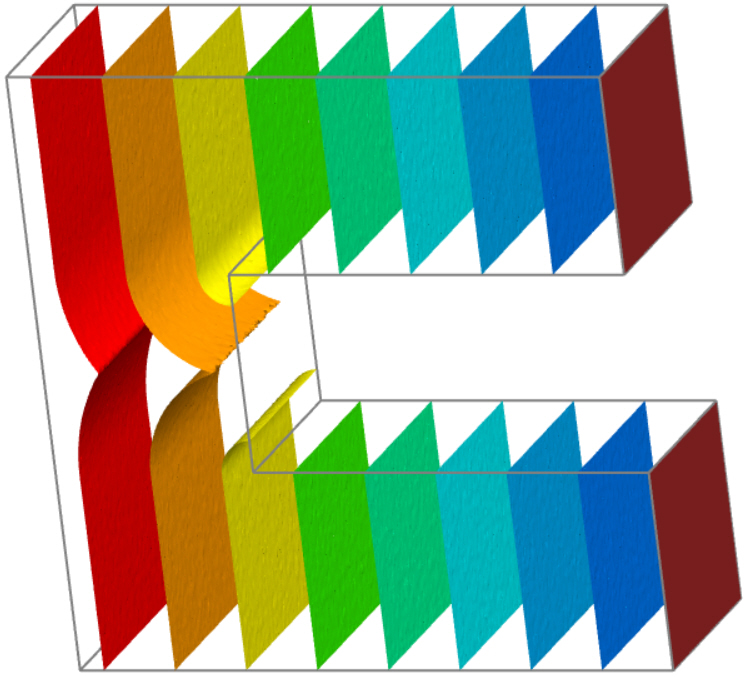} \\
			\ref{ex_3} & \ref{ex_4} \\
			\includegraphics[height=3cm]{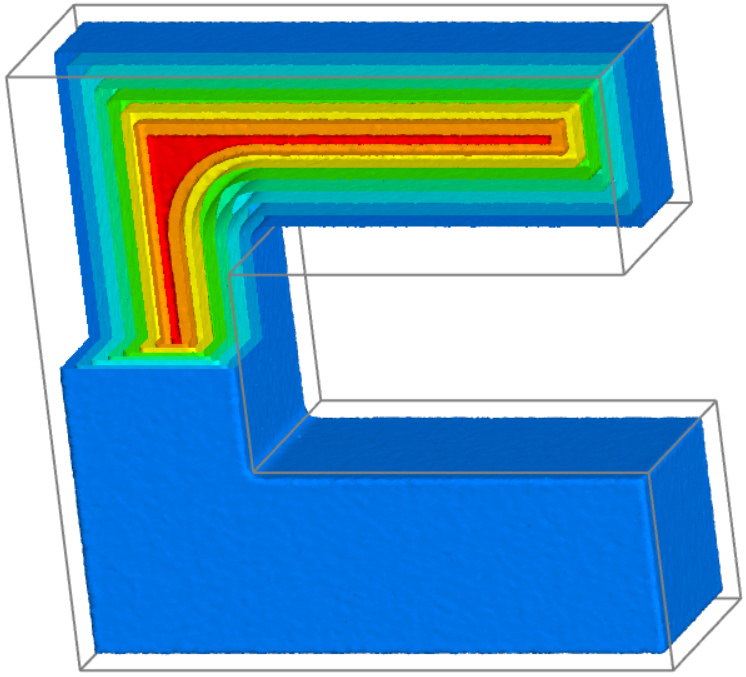} &
			\includegraphics[height=3cm]{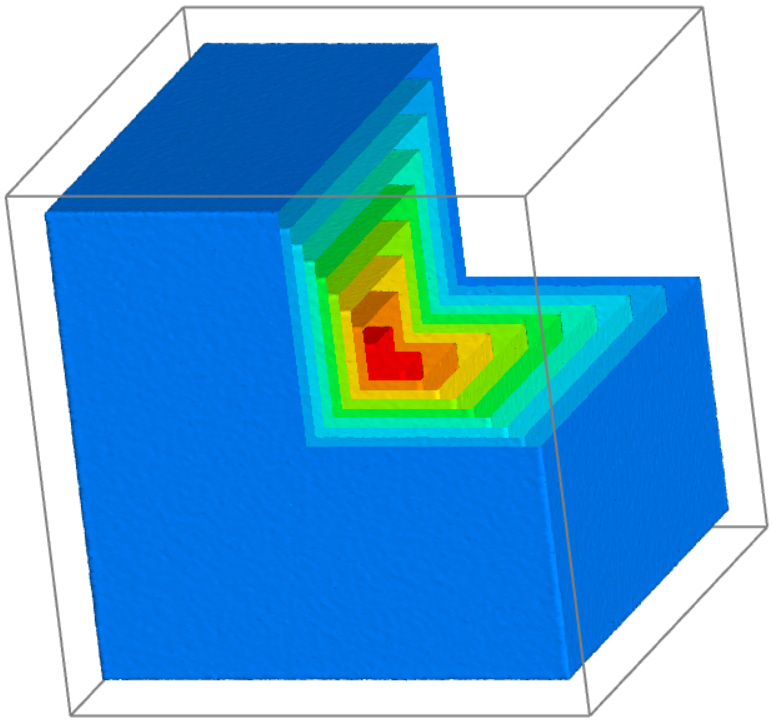} \\
			\ref{ex_5} & \ref{ex_6} \\
			\includegraphics[height=3cm]{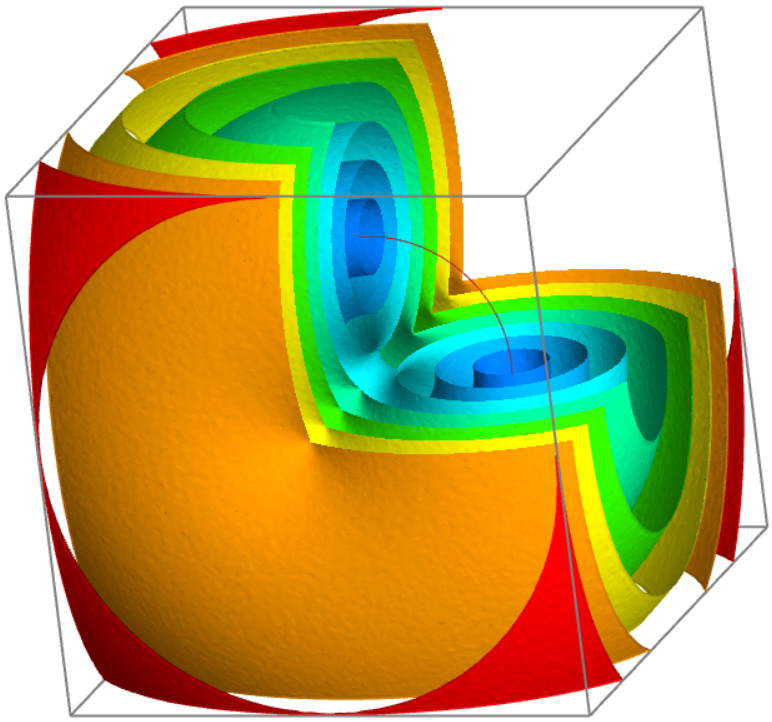} &
			\includegraphics[height=3cm]{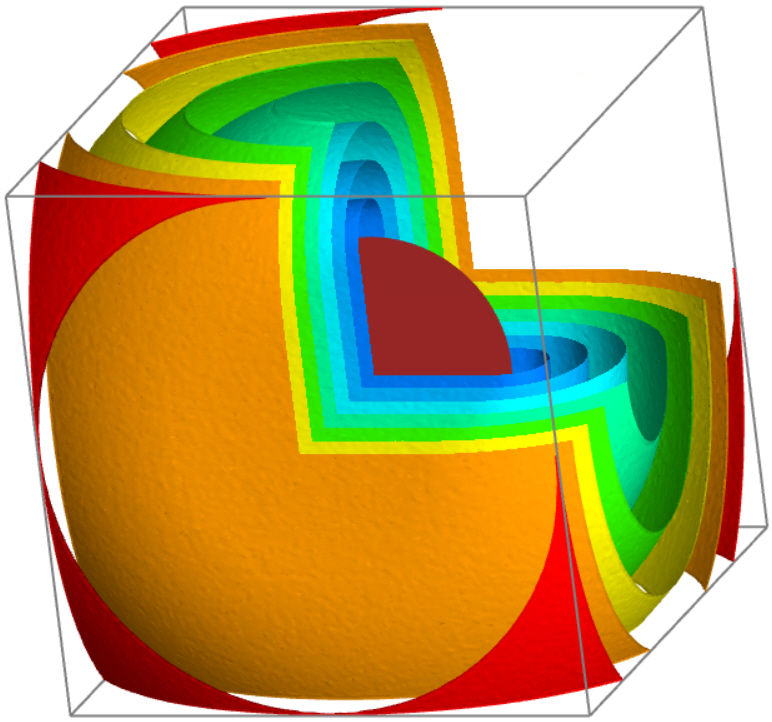} \\
			\ref{ex_7} & \ref{ex_8} \\
			\includegraphics[height=3cm]{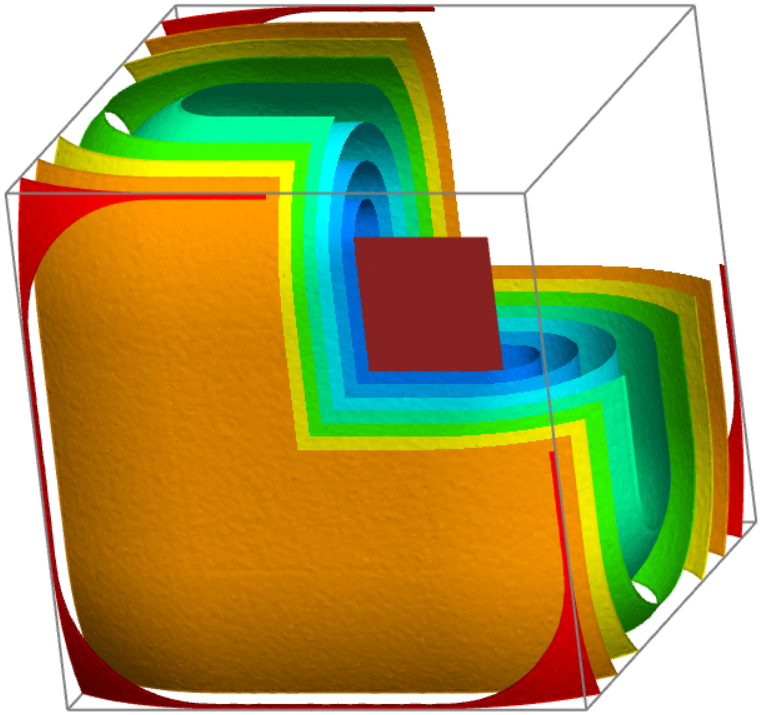} &
			\includegraphics[height=3cm]{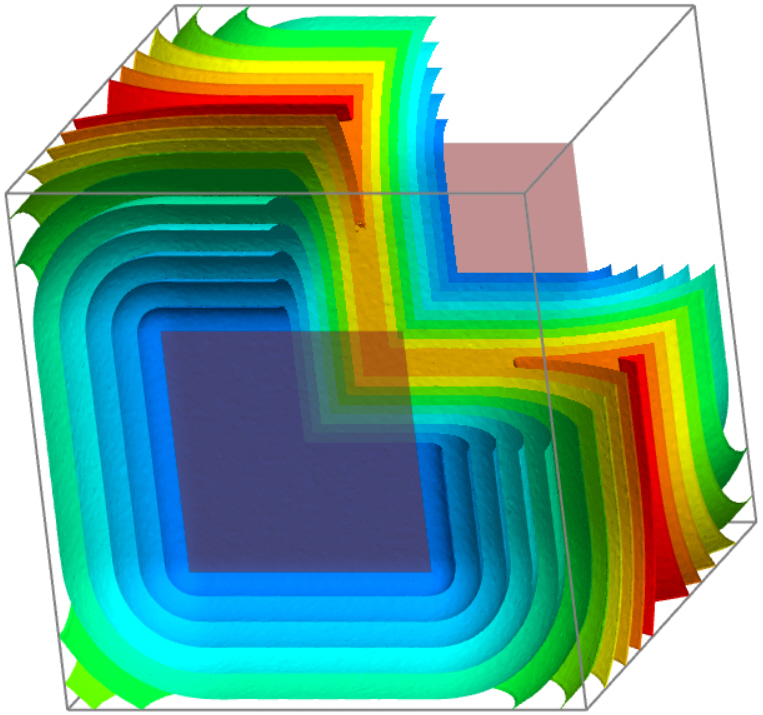} \\
			\ref{ex_9} & \ref{ex_10} \\
		\end{tabular}
	\end{center}
	\caption{Iso-surfaces of numerical solutions computed by the proposed algorithm~\eqref{eq:prop_alg} are presented on the level $L=4$ in Table~\ref{tab:meshes}.} \label{fig:isosurf}
\end{figure}

Prior to the numerical properties of the proposed algorithm, equidistant isosurfaces of numerical solutions computed by the proposed algorithm~\eqref{eq:prop_alg} are presented in Figure~\ref{fig:isosurf} on the level $L=4$ in Table~\ref{tab:meshes}. They are qualitatively shown as a distance function from a given $\Gamma$ illustrated by the color of dark red. In the cases of~\ref{ex_1},~\ref{ex_2}, and~\ref{ex_10}, we use a transparency on $\Gamma$ to visually observe isosurfaces behind~$\Gamma$. In the cases of~\ref{ex_5} and~\ref{ex_6}, the surface $\Gamma$ is not presented because $\Gamma = \partial \Omega$.

\begin{table}
	\centering
	\small
	\begin{tabular}{cS[table-format=1]S[table-format=1.2,
			table-figures-exponent=2,
			table-sign-mantissa,
			table-sign-exponent]S[table-format=1.2]S[table-format=1.2,
			table-figures-exponent=2,
			table-sign-mantissa,
			table-sign-exponent]S[table-format=1.2]}
		\toprule
		\multicolumn{1}{c}{} & \multicolumn{1}{c}{L} & \multicolumn{1}{c}{$E^1$}    & \multicolumn{1}{c}{$EOC$} & \multicolumn{1}{c}{$E^{\infty}$} & \multicolumn{1}{c}{$EOC$} \\
		\midrule
		\multirow{4}{*}{\ref{ex_1}} & 1 & 2.00E-03 & 2.20  & 6.14E-03     & 0.80  \\
		& 2 & 5.79E-04 & 1.29  & 3.90E-03     & 0.96  \\
		& 3 & 2.66E-04 & 2.29  & 2.19E-03     & 0.91  \\
		& 4 & 6.19E-05 &       & 1.22E-03     &       \\
		\cmidrule(lr){1-2} \cmidrule(lr){3-4} \cmidrule(lr){5-6}
		\multirow{4}{*}{\ref{ex_2}} & 1 & 2.76E-03 & 1.27  & 3.72E-02     & 0.71  \\
		& 2 & 1.39E-03 & 1.69  & 2.53E-02     & 1.15  \\
		& 3 & 5.43E-04 & 1.59  & 1.33E-02     & 1.29  \\
		& 4 & 2.11E-04 &       & 6.17E-03     &       \\
		\cmidrule(lr){1-2} \cmidrule(lr){3-4} \cmidrule(lr){5-6}
		\multirow{4}{*}{\ref{ex_3}} & 1 & 1.28E-02 & 1.54  & 3.68E-02     & 1.28  \\
		& 2 & 5.85E-03 & 1.09  & 1.91E-02     & 1.16  \\
		& 3 & 3.21E-03 & 1.31  & 1.01E-02     & 1.08  \\
		& 4 & 1.63E-03 &       & 5.79E-03     &       \\
		\cmidrule(lr){1-2} \cmidrule(lr){3-4} \cmidrule(lr){5-6}
		\multirow{4}{*}{\ref{ex_4}} & 1 & 3.12E-03 & 1.24  & 5.54E-02     & 1.03  \\
		& 2 & 1.66E-03 & 1.07  & 3.28E-02     & 1.01  \\
		& 3 & 9.19E-04 & 1.26  & 1.89E-02     & 1.40  \\
		& 4 & 4.76E-04 &       & 9.12E-03     &       \\
		\cmidrule(lr){1-2} \cmidrule(lr){3-4} \cmidrule(lr){5-6}
		\multirow{4}{*}{\ref{ex_5}} & 1 & 5.88E-03 & 1.04  & 5.72E-02     & 0.90  \\
		& 2 & 3.46E-03 & 1.89  & 3.62E-02     & 1.12  \\
		& 3 & 1.23E-03 & 2.24  & 1.96E-02     & 0.92  \\
		& 4 & 3.84E-04 &       & 1.21E-02     &       \\ 
		\cmidrule(lr){1-2} \cmidrule(lr){3-4} \cmidrule(lr){5-6}		
		\multirow{4}{*}{\ref{ex_6}} & 1 & 6.90E-02 & 3.44 & 6.20E-01 & 2.76 \\
		& 2 & 7.53E-03 & 2.05 & 1.05E-01 & 1.72 \\
		& 3 & 1.87E-03 & 1.66 & 3.27E-02 & 0.55 \\
		& 4 & 5.98E-04 & ~ & 2.24E-02 &  \\
		\cmidrule(lr){1-2} \cmidrule(lr){3-4} \cmidrule(lr){5-6}
		\multirow{4}{*}{\ref{ex_7}} & 1 & 3.12E-01 & 1.41 & 6.66E-01 & 1.63 \\
		& 2 & 1.26E-01 & 1.97 & 2.33E-01 & 2.16 \\
		& 3 & 3.32E-02 & 2.27 & 5.38E-02 & 1.76 \\
		& 4 & 6.93E-03 & ~ & 1.60E-02 & \\ 
		\cmidrule(lr){1-2} \cmidrule(lr){3-4} \cmidrule(lr){5-6}
		\multirow{4}{*}{\ref{ex_8}} & 1 & 2.89E-01 & 1.48 & 6.53E-01 & 1.62 \\
		& 2 & 1.11E-01 & 2.00 & 2.30E-01 & 2.15 \\
		& 3 & 2.86E-02 & 2.31 & 5.35E-02 & 2.40 \\
		& 4 & 5.84E-03 &  & 1.03E-02 \\
		\cmidrule(lr){1-2} \cmidrule(lr){3-4} \cmidrule(lr){5-6}
		\multirow{4}{*}{\ref{ex_9}}  & 1 & 2.60E-01 & 1.40  & 7.81E-01     & 1.39  \\
		& 2 & 1.06E-01 & 1.84  & 3.20E-01     & 1.98  \\
		& 3 & 3.03E-02 & 2.16  & 8.33E-02     & 2.02  \\
		& 4 & 6.83E-03 &       & 2.07E-02     &       \\
		\cmidrule(lr){1-2} \cmidrule(lr){3-4} \cmidrule(lr){5-6}
		\multirow{4}{*}{\ref{ex_10}} & 1 & 1.01E+00 & 1.88  & 2.06E+00     & 1.60  \\
		& 2 & 3.03E-01 & 1.74  & 7.35E-01     & 1.63  \\
		& 3 & 9.26E-02 & 1.41  & 2.43E-01     & 1.35  \\
		& 4 & 3.50E-02 &       & 9.64E-02     & \\
		\bottomrule 
	\end{tabular}
	\caption{The $EOC$s~\eqref{eq:EOC} of all examples for a numerical solution of~\eqref{eq:prop_alg} with $\epsilon_n$, $n=5$, are presented. $\text{L}$ is the level of mesh listed in Table~\ref{tab:meshes}.}\label{tab:EOC}
\end{table}

The first numerical property is an experimental order of convergence ($EOC$). Since exact solutions for all examples are known, we compute the errors $E^1_{\text{L}}$ and $E^{\infty}_{\text{L}}$ of $L^{1}$ and $L^{\infty}$ norms between a numerical solution on the $\text{L}^\text{th}$ level of mesh and an exact solution, respectively. Then, for each error, the corresponding $EOC$ is calculated by
\begin{align}\label{eq:EOC}
	EOC_{\text{L}} = \frac{\log\left(\frac{E_\text{L+1}}{E_\text{L}}\right)}{\log\left(\frac{h_\text{L+1}}{h_\text{L}}\right)}, \quad \text{L} \in \{1,\:2,\:3\}.
\end{align}
In Table~\ref{tab:EOC}, we present $EOC$s of all examples for a numerical solution of the proposed algorithm~\eqref{eq:prop_alg} with $\epsilon_n$, $n=5$. For smooth solutions of~\ref{ex_8} and ~\ref{ex_9}, the $EOC$s with $E^1$ and $E^{\infty}$ errors are larger than $2$. In~\ref{ex_1}, the $EOC$s with $E^1$ is larger than $2$, but the $EOC$s with $E^{\infty}$ is close to $1$ because of a singularity at the origin. For all non-smooth solutions, the $EOC$s with $E^1$ and $E^{\infty}$ errors are close to $1$. Compared to the $EOC$s in~\cite{ref:HMFB21}, the behavior of $EOC$ is quite similar.

\begin{table}
	\centering
	\small
	\begin{tabular}{S[table-format=1]cS[table-format=1.2,
			table-figures-exponent=2,
			table-sign-mantissa,
			table-sign-exponent]S[table-format=1.2]S[table-format=1.2,
			table-figures-exponent=2,
			table-sign-mantissa,
			table-sign-exponent]S[table-format=1.2]}
		\toprule
		\multicolumn{1}{c}{L} & \multicolumn{1}{c}{$\epsilon_n$} & \multicolumn{1}{c}{$E^1$}    & \multicolumn{1}{c}{$EOC$} & \multicolumn{1}{c}{$E^{\infty}$} & \multicolumn{1}{c}{$EOC$} \\
		\midrule		
		\cellcolor{gray4} 1 & \multirow{4}{*}{$h_{\text{L}}^1$} & \cellcolor{gray4} 9.06E-02 & 0.33  & \cellcolor{gray4}  2.68E-01     & 0.34  \\
		\cellcolor{gray3} 2 & & \cellcolor{gray3} 7.50E-02 & 0.75  & \cellcolor{gray3} 2.21E-01     & 0.66  \\
		\cellcolor{gray2} 3 & & \cellcolor{gray2} 4.76E-02 & 0.80  & \cellcolor{gray2} 1.48E-01     & 0.55  \\
		\cellcolor{gray1} 4 & & \cellcolor{gray1} 2.86E-02 &       & \cellcolor{gray1} 1.04E-01     &       \\
		\cmidrule(lr){1-2} \cmidrule(lr){3-4} \cmidrule(lr){5-6}
		\cellcolor{gray4} 1 & \multirow{4}{*}{$h_{\text{L}}^\frac{3}{2}$} & \cellcolor{gray4} 2.11E-02 & 0.82  & \cellcolor{gray4} 5.57E-02     & 0.96  \\
		\cellcolor{gray3} 2 & & \cellcolor{gray3} 1.33E-02 & 1.19  & \cellcolor{gray3} 3.23E-02     & 0.71  \\
		\cellcolor{gray2} 3 & & \cellcolor{gray2} 6.48E-03 & 1.14  & \cellcolor{gray2} 2.10E-02     & 0.40  \\
		\cellcolor{gray1} 4 & & \cellcolor{gray1} 3.14E-03 &       & \cellcolor{gray1} 1.63E-02     &       \\
		\cmidrule(lr){1-2} \cmidrule(lr){3-4} \cmidrule(lr){5-6}
		\cellcolor{gray4} 1 & \multirow{4}{*}{$h_{\text{L}}^2$} & \cellcolor{gray4} 6.81E-03 & 1.54  & \cellcolor{gray4} 1.80E-02     & 0.95  \\
		\cellcolor{gray3} 2 & & \cellcolor{gray3} 2.86E-03 & 2.09  & \cellcolor{gray3} 1.05E-02     & 2.59  \\
		\cellcolor{gray2} 3 & & \cellcolor{gray2} 8.08E-04 & 1.38  & \cellcolor{gray2} 2.20E-03     & 0.09  \\
		\cellcolor{gray1} 4 & & \cellcolor{gray1} 3.35E-04 &       & \cellcolor{gray1} 2.07E-03     &       \\
		\cmidrule(lr){1-2} \cmidrule(lr){3-4} \cmidrule(lr){5-6}
		\cellcolor{gray4} 1 & \multirow{4}{*}{$h_{\text{L}}^\frac{5}{2}$} & \cellcolor{gray4} 2.00E-03 & 2.20  & \cellcolor{gray4} 6.14E-03     & 0.80  \\
		\cellcolor{gray3} 2 & & \cellcolor{gray3} 5.79E-04 & 1.29  & \cellcolor{gray3} 3.90E-03     & 0.96  \\
		\cellcolor{gray2} 3 & & \cellcolor{gray2} 2.66E-04 & 2.29  & \cellcolor{gray2} 2.19E-03     & 0.91  \\
		\cellcolor{gray1} 4 & & \cellcolor{gray1} 6.19E-05 &       & \cellcolor{gray1} 1.22E-03     &      \\
		\bottomrule 
	\end{tabular}
	\caption{For the case of~\ref{ex_1}, errors $E^1$ and $E^{\infty}$ of numerical solutions~\eqref{eq:prop_alg} with $\epsilon_n$ from $n=2$ to $n=5$ are presented on all levels of meshes. From a fixed $\epsilon_n$, the $EOC$s are also shown on different levels of meshes.}\label{tab:sp_full}
\end{table}

The second numerical property is the behavior of the errors versus the regularization parameter $\epsilon_n$ on a fixed level of meshes. For each $n$ on the $\text{L}^{\text{th}}$ level of mesh, the proposed algorithm~\eqref{eq:prop_alg} provides a numerical solution $u^{n}$ with $\epsilon_n = (h_\text{L})^{\frac{1}{2}n}$. For the next $n+1$, we use the solution $u^n$ and then find the next solution $u^{n+1}$ with $\epsilon_{n+1}$ ($< \epsilon_n$). In Table~\ref{tab:sp_full}, for the case of~\ref{ex_1}, errors $E^1$ and $E^{\infty}$ of numerical solutions $u^n$ with $\epsilon_n$ from $n=2$ to $n=5$ are presented on all levels of meshes. A crucial observation is that the choice of $\epsilon_5 = h_L^{\frac{5}{2}}$ brings a better result, that is, smaller errors, than the other regularization values $\epsilon_n$ for $1 \leq n \leq 4$ or $n=6$. Since we use $K_1=1$ in~\eqref{eq:prop_alg}, the results of $n=1$ are far from the exact solution. On a fixed level of mesh, when the regularization parameter $\epsilon_n$ is smaller, that is, $n$ becomes larger, the errors $E^1$ and $E^{\infty}$ become smaller until $n=5$. The mentioned property can be seen on the rows with the same gray color in Table~\ref{tab:sp_full}. For example, when $\text{L}=1$, by the value on the second row of $E^1$ column, the error on every fourth row below in the same column decreases; see the error values shadowed by the darkest gray color in Table~\ref{tab:sp_full}. Also, the $EOC$s on different levels of meshes become better from $n=2$ to $n=5$. When $n \geq 6$, the effect of the Laplacian regularizer is too small to solve the linear system~\eqref{eq:prop_alg} stably enough. A similar instability of using too small regularization parameter is also observed in~\cite{ref:THP02,ref:CWW13,ref:CWW17}. A relation between the regularization parameter and the order of numerical scheme is also observed in~\cite{ref:GR09}. A further numerical analysis is necessary to find an optimal regularization parameter to minimize an error between a numerical solution on a discrete space of~\eqref{eq:EkEq_Lap} and a viscosity solution of~\eqref{eq:EkEq_org}, which is out of the scope of this paper. 

The third numerical property is a comparison of computational cost. To minimize a systematical bias, we purposely choose the time-relaxed bidirectional eikonal equation~\cite{ref:HMFB21} already implemented in \texttt{AVL FIRE}\textsuperscript{\texttt{TM}}. The proposed algorithm is also implemented by the same language (Fortran 2003) and all algorithms are compiled by the same compiler options. 

\begin{table}
	\centering
	\small
	\begin{tabular}{c
			S[table-format=1]
			S[table-format=4.2]
			S[table-format=2]
			S[table-format=1.3]
			S[table-format=5.2]
			S[table-format=4]
			S[table-format=2.3]}
		\toprule
		& 
		\multicolumn{1}{c}{L} & 
		\multicolumn{1}{c}{Time$^{\text{1}}$(s)} & 
		\multicolumn{1}{c}{$\displaystyle \sum_{n=1}^5 K_n$} &  
		\multicolumn{1}{c}{Final $T$} & \multicolumn{1}{c}{Time$^{\text{2}}$(s)}  & 
		\multicolumn{1}{c}{$N_{tot}$} & 
		\multicolumn{1}{c}{Ratio}  \\
		\cmidrule(lr){1-2}\cmidrule(lr){3-4}\cmidrule(lr){5-7}\cmidrule(lr){8-8}
		\multirow{4}{*}{\ref{ex_1}} & 1 & 24.32   & 36 & 1.400 & 52.10    & 35  & 2.142  \\
		& 2 & 60.01   & 27 & 1.400 & 328.67   & 70  & 5.477  \\
		& 3 & 170.64  & 20 & 1.390 & 2253.86  & 139 & 13.208 \\
		& 4 & 346.83  & 8  & 1.505 & 17280.66 & 301 & 49.825 \\
		\midrule  
		\multirow{4}{*}{\ref{ex_2}} & 1 & 88.09   & 33 & 1.880 & 291.08   & 47  & 3.304  \\
		& 2 & 171.21  & 21 & 1.740 & 1586.78  & 87  & 9.268  \\
		& 3 & 351.49  & 11 & 1.770 & 10332.52 & 177 & 29.396 \\
		& 4 & 1049.69 & 6  & 1.845 & 72433.44 & 369 & 69.005 \\		
		\midrule
		\multirow{4}{*}{\ref{ex_3}} & 1 & 18.14   & 45 & 4.120 & 55.31   & 103 & 3.049  \\
		& 2 & 30.67   & 28 & 4.260 & 247.40  & 213 & 8.067  \\
		& 3 & 51.41   & 14 & 4.200 & 1300.40 & 420 & 25.295 \\
		& 4 & 139.66  & 10 & 4.225 & 7609.97 & 845 & 54.488 \\
		\midrule  
		\multirow{4}{*}{\ref{ex_4}} & 1 & 14.26   & 30 & 2.680 & 36.86   & 67  & 2.584  \\
		& 2 & 23.40   & 19 & 3.160 & 186.60  & 158 & 7.973  \\
		& 3 & 42.07   & 10 & 2.930 & 921.51  & 293 & 21.906 \\
		& 4 & 119.39  & 7  & 2.875 & 5183.78 & 575 & 43.419 \\
		\midrule  
		\multirow{4}{*}{\ref{ex_5}} & 1 & 10.92   & 20 & \multicolumn{1}{c}{$T_{M}=2$} & 27.76    & 50  & 2.543  \\
		& 2 & 20.47   & 16 & \multicolumn{1}{c}{$T_{M}=2$} & 119.38   & 100 & 5.832  \\
		& 3 & 46.97   & 13 & 0.770                         & 245.93   & 77  & 5.236  \\
		& 4 & 124.66  & 9  & 0.625                         & 1185.81  & 125 & 9.513 \\
		\bottomrule
	\end{tabular}
	\caption{For all examples, a comparison of computational cost is presented by using $2^{\text{L+2}}$ numbers of CPUs on the $\text{L}^{\text{th}}$ level of mesh. Time$^\text{1}$ and Time$^\text{2}$ are the computation time in seconds of the proposed algorithm~\eqref{eq:prop_alg} and the algorithm in~\cite{ref:HMFB21}, respectively, and the corresponding total number of iterations are shown right next to the computational time. The final $T$ to solve~\eqref{eq:EkEq_Time} is decided by the same $E^1$ error value as the proposed method; see more details in Section~\ref{sec:num}.}\label{tab:time_cost}
\end{table}

Since the time-relaxed bidirectional eikonal equation is time-dependent and the governing equation in this paper is time-independent, we stop the time evolution in~\eqref{eq:EkEq_Time} right before the $E^1$ error of~\eqref{eq:EkEq_Time} becomes smaller than the $E^1$ error of the proposed algorithm. That is, we measure a computational cost until two methods reach the same error bound. In Table~\ref{tab:time_cost}, such a final time $T$ is shown on the column labeled by ``Final $T$'' for all examples. On that column, $T_M$ means that $E^1$ error of~\eqref{eq:EkEq_Time} is not smaller than the $E^1$ error of the proposed algorithm until the predetermined final time $T_M$, specified in~\cite{ref:HMFB21}. Time$^\text{1}$ and Time$^\text{2}$ are the computation time in seconds for the proposed algorithm~\eqref{eq:prop_alg} and the algorithm in~\cite{ref:HMFB21}, respectively, and the corresponding total number of iterations are shown right next to the computational time. The calculations of using $2^\text{L+2}$ numbers of CPUs for all examples in the $\text{L}^{\text{th}}$ level of mesh are repeated five times in a cluster, a distributed system (Intel Xeon$^{\text{\textregistered}}$ Gold 6154 Processor $3.00$GHz $20$ CPUs and $20$ gigabyte memory) and the computational time (Time$^\text{1}$ and Time$^\text{2}$) in Table~\ref{tab:time_cost} is the average of five measurements. Since the distance information in~\eqref{eq:EkEq_Time} is evolved from $\Gamma$ over time, the time-relaxed bidirectional equation has certainly a disadvantage in computational time whenever it is necessary to compute a distance further away from $\Gamma$. On the last column, it shows how much the proposed algorithm is faster than the algorithm to solve the time-relaxed bidirectional eikonal equation to reach the same $E^1$ error. A crucial point of the last column is that the computational efficiency in the proposed algorithm becomes better when there are more number of cells in a mesh.

In the case of~\ref{ex_5}, the time ratio is quite different by other examples because $\Gamma = \partial \Omega$ makes the traveling distance much shorter than other examples. In other words, the computational time of the proposed algorithm becomes faster than the previous approach~\cite{ref:HMFB21} as long as the region of interest to find distance values are far away from $\Gamma$.

\section{Conclusion}\label{sec:conclusion}
We present a cell-centered finite volume method to solve a Laplacian regularized eikonal equation with Soner boundary condition on polyhedral meshes in order to compute a distance function from given objects. Using a linearized form of the equation, a numerical solution is sequentially updated by a decreasing sequence of the regularization parameters depending on a characteristic length of discretized domain. The normalized gradient field of the first solution in the sequence is substantially improved on most part of domain. As the characteristic length becomes smaller, the regularization parameter becomes smaller and a convergence to the viscosity solution is numerically verified. The $EOC$ of $L^1$ norm of the error is shown to be the second order for tested smooth solutions. Compared to the computational time of solving the time-relaxed bidirectional eikonal equation, the proposed algorithm has an advantage to dramatically reduce the time when a larger number of cells is used or a region of interest is far away from where the distance measurement starts. The implementation of parallel computing using domain decomposition with the $1$-ring face neighbor structure can be done straightforwardly by a standard cell-centered finite volume code. 

\section*{Acknowledgments}
The authors thank Prof. Silvia Tozza in University of Bologna, Italy, for comments of Soner boundary condition. We also sincerely thank Dr. Branislav Basara and Dr. Reinhard Tatschl in AVL List GmbH, Austria, for supporting the University Partnership Program\footnote{See more details in AVL Advanced Simulation Technologies University Partnership Program: \href{https://www.avl.com/documents/10138/3372587/AVL_UPP_Flyer.pdf}{https://www.avl.com/documents/10138/3372587/AVL\_UPP\_Flyer.pdf}}.

%\section*{References}

%\bibliographystyle{splncs}
\bibliographystyle{plain}
\bibliography{reference}

\begin{thebibliography}{10}

\bibitem{ref:AA14}
G.~Aubert and J.-F. Aujol.
\newblock Poisson skeleton revisited: a {N}ew mathematical perspective.
\newblock {\em Journal of Mathematical Imaging and Vision}, 48:149--159, 2014.

\bibitem{ref:BB91}
B.~Baldwin and T.~Barth.
\newblock A one-equation turbulence transport model for high reynolds number
  wall-bounded flows.
\newblock American Institute of Aeronautics and Astronautics 29th Aerospace
  Sciences Meetin, 91-0610, 1991.

\bibitem{ref:BL78}
B.~Baldwin and H.~Lomax.
\newblock Thin-layer approximation and algebraic model for separated
  turbulentflows.
\newblock American Institute of Aeronautics and Astronautics 16th Aerospace
  Sciences Meeting, 78-257, 1978.

\bibitem{ref:BS98}
T.~J. Barth and J.~A. Sethian.
\newblock Numerical schemes for the {H}amilton-{J}acobi and level set equations
  on triangulated domains.
\newblock {\em Journal of Computational Physics}, 145:1--40, 1998.

\bibitem{ref:BF15}
A.~G. Belyaev and P.-A. Fayolle.
\newblock On variational and {PDE}-based distance function approximations.
\newblock {\em Computer Graphics Forum}, 34(8):104--118, 2015.

\bibitem{ref:BF19}
A.~G. Belyaev and P.-A. Fayolle.
\newblock A variational method for accurate distance function estimation.
\newblock pages 175--181. Numerical Geometry, Grid Generation and Scientific
  Computing, Springer International Publishing, 2019.

\bibitem{ref:BHS84}
K.~B{\"o}hmer, P.~W. Hemker, and H.~J. Stetter.
\newblock The defect correction approach.
\newblock In {\em Defect correction methods}, pages 1--32. Springer, 1984.

\bibitem{ref:CGP15}
A.~Caboussat, R.~Glowinski, and T.-W. Pan.
\newblock On the numerical solution of some eikonal equations: An elliptic
  solver approach.
\newblock {\em Chin. Ann. Math. Ser. B}, 36:689--702, 2015.

\bibitem{ref:CDL90}
I.~Capuzzo-Dolcetta and P.-L. Lions.
\newblock Hamilton-{J}acobi equations with state constraints.
\newblock {\em Transactions of the American Mathematical Society}, 318:1990,
  643-683.

\bibitem{ref:CV19}
A.~G. Churbanov and P.~N. Vabishchevich.
\newblock Numerical solution of boundary value problems for the eikonal
  equation in an anisotropic medium.
\newblock {\em Journal of Computational and Applied Mathematics}, 362:55--67,
  2019.

\bibitem{ref:CFG93}
P.~Colli-Franzone and L.~Guerri.
\newblock Spreading of excitation in 3-{D} models of the anisotropic cardiac
  tissue. {I}. validation of the eikonal model.
\newblock {\em Mathematical Biosciences}, 113:145--209, 1993.

\bibitem{ref:CEL84}
M.~G. Crandall, L.~C. Evans, and P.-L. Lions.
\newblock Some properties of viscosity solutions of {H}amilton-{J}acobi
  equations.
\newblock {\em Transactions of the American Mathematical Society}, 282:1984,
  487-502.

\bibitem{ref:CWW13}
K.~Crane, C.~Weischedel, and M.~Wardetzky.
\newblock Geodesics in heat: A new approach to computing distance based on heat
  flow.
\newblock {\em ACM Transactions on Graphics (TOG)}, 32:152:1--152:11, 2013.

\bibitem{ref:CWW17}
K.~Crane, C.~Weischedel, and M.~Wardetzky.
\newblock The heat method for distance computation.
\newblock {\em Communications of the ACM}, 60:90--99, 2017.

\bibitem{ref:DES11}
K.~Deckelnick, C.~M. Elliott, and V.~Styles.
\newblock Numerical analysis of an inverse problem for the eikonal equation.
\newblock {\em Numerische Mathematik}, 119:245--269, 2011.

\bibitem{ref:EIN21}
H.~Ennaji, N.~Igbida, and V.~T. Nguyen.
\newblock Augmented lagrangian methods for degenerate {H}amilton-{J}acobi
  equations.
\newblock {\em Calculus of Variations and Partial Differential Equations},
  60:238, 2021.

\bibitem{ref:E98}
L.~C. Evans.
\newblock {\em Partial differential equations}.
\newblock American Methematical Society, Providence, R.I., 1998.

\bibitem{ref:FT09}
M.~Falcone and C.~Truini.
\newblock A level-set algorithm for front propagation in the presence of
  obstacles.
\newblock {\em Rendiconti di Matematica e delle sue Applicazioni}, 29:29--50,
  2009.

\bibitem{ref:FS02}
E.~Fares and W.~Schr\"oder.
\newblock A differential equation for approximate wall distance.
\newblock {\em International Journal for Numerical Methods in Fluids},
  39:743--762, 2002.

\bibitem{ref:FB20}
P.-A. Fayolle and A.~G. Belyaev.
\newblock An {ADMM}-based scheme for distance function approximation.
\newblock {\em Numerical Algorithms}, 84:983--996, 2020.

\bibitem{ref:FMHMB19}
P.~Frolkovi{\v{c}}, K.~Mikula, J.~Hahn, D.~Martin, and B.~Basara.
\newblock Flux balanced approximation with least-squares gradient for diffusion
  equation on polyhedral mesh.
\newblock {\em Discrete \& Continuous Dynamical Systems - S}, 2020.

\bibitem{ref:FJPKW11}
Z.~Fu, W.-K. Jeong, Y.~Pan, R.~M. Kirby, and R.~T. Whitaker.
\newblock A fast iterative method for solving the eikonal equation on
  triangulated surfaces.
\newblock {\em SIAM Journal on Scientific Computing}, 33:2468--2488, 2011.

\bibitem{ref:FKW13}
Z.~Fu, R.~M. Kirby, and R.~T. Whitaker.
\newblock A fast iterative method for solving the eikonal equation on
  tetrahedral domains.
\newblock {\em SIAM Journal on Scientific Computing}, 35:C473--C494, 2013.

\bibitem{ref:GK06}
P.~A. Gremaud and C.~M. Kuster.
\newblock Computational study of fast methods for the eikonal equation.
\newblock {\em SIAM Journal on Scientific Computing}, 27:1803–1816, 2006.

\bibitem{ref:GR09}
K.~S. Gurumoorthy and A.~Rangarajan.
\newblock A {S}chr\"{o}dinger equation for the fast computation of approximate
  euclidean distance functions.
\newblock pages 100--111. Scale Space and Variational Methods in Computer
  Vision. SSVM 2009. Lecture Notes in Computer Science, vol 5567, Springer,
  Berlin, Heidelberg, 2009.

\bibitem{ref:HMFB17}
J.~Hahn, K.~Mikula, P.~Frolkovi{\v{c}}, and B.~Basara.
\newblock Inflow-based gradient finite volume method for a propagation in a
  normal direction in a polyhedron mesh.
\newblock {\em Journal of Scientific Computing}, 72:442--465, 2017.

\bibitem{ref:HMFB17_1}
J.~Hahn, K.~Mikula, P.~Frolkovi{\v{c}}, and B.~Basara.
\newblock Semi-implicit level set method with inflow-based gradient in a
  polyhedron mesh.
\newblock In C.~Canc{\`e}s and P.~Omnes, editors, {\em Finite Volumes for
  Complex Applications VIII - Hyperbolic, Elliptic and Parabolic Problems},
  pages 81--89. Springer International Publishing, 2017.

\bibitem{ref:HMFMB19}
J.~Hahn, K.~Mikula, P.~Frolkovi{\v{c}}, M.~Medl\char39a, and B.~Basara.
\newblock Iterative inflow-implicit outflow-explicit finite volume scheme for
  level-set equations on polyhedron meshes.
\newblock {\em Computers \& Mathematics with Applications}, 77:1639--1654,
  2019.

\bibitem{ref:HMFB21}
Jooyoung Hahn, Karol Mikula, Peter Frolkovi\v{c}, and Branislav Basara.
\newblock Finite volume method with the {S}oner boundary condition for
  computing the signed distance function on polyhedral meshes.
\newblock {\em International Journal for Numerical Methods in Engineering},
  123:1057--1077, 2022.

\bibitem{ref:HT05}
S.-R. Hysing and S.~Turek.
\newblock The eikonal equation: {N}umerical efficiency vs. algorithmic
  complextiy on quadrilateral grids.
\newblock pages 22--31. Proceedings of ALGORITMY, 2005.

\bibitem{ref:JW08}
W.-K. Jeong and R.~T. Whitaker.
\newblock A fast iterative method for eikonal equations.
\newblock {\em SIAM Journal on Scientific Computing}, 30:2512--2534, 2008.

\bibitem{ref:K91}
J.~P. Keener.
\newblock An eikonal-curvature equation for action potential propagation in
  myocardium.
\newblock {\em Journal of Mathematical Biology}, 29:629--651, 1991.

\bibitem{ref:KS98}
R.~Kimmel and J.~A. Sethian.
\newblock Computing geodesic paths on manifolds.
\newblock {\em Proceedings of the National Academy of Sciences}, 95:8431--8435,
  1998.

\bibitem{ref:M16}
A.~Manz.
\newblock {\em Modeling of End-Gas Autoignition for Knock Prediction in
  Gasoline Engines}.
\newblock Logos Verlag Berlin, 2016.

\bibitem{ref:OF00}
S.~Osher and R.~Fedkiw.
\newblock {\em Level set methods and dynamic implicit surfaces}.
\newblock Springer, Berlin, 2000.

\bibitem{ref:P04}
M.~Peri{\v c}.
\newblock Flow simulation using control volumens of arbitrary polyhedral shape.
\newblock pages 25--29. ERCOFTAC Bulletin 62, 2004.

\bibitem{ref:P00}
N.~Peters.
\newblock {\em Turbulent Combustion}.
\newblock Cambridge Monographs on Mechanics. Cambridge University Press, 2000.

\bibitem{ref:PAS95}
M.~A. Price, C.~G. Armstrong, and M.~A. Sabin.
\newblock Hexahedral mesh generation by medial surface subdivision: Part i.
  solids with convex edges.
\newblock {\em International Journal for Numerical Methods in Engineering},
  38:3335--3359, 1995.

\bibitem{ref:QZZ07a}
J.-L. Qian, Y.-T. Zhang, and H.-K. Zhao.
\newblock Fast sweeping methods for eikonal equations on triangular meshes.
\newblock {\em SIAM Journal on Numerical Analysis}, 31:83--107, 2007.

\bibitem{ref:QRPG04}
W.~R. Quadros, K.~Ramaswami, F.~B. Prinz, and B.~Gurumoorthy.
\newblock Laytracks: a new approach to automated geometry adaptive
  quadrilateral mesh generation using medial axis transform.
\newblock {\em International Journal for Numerical Methods in Engineering},
  61:209--237, 2004.

\bibitem{ref:RS05}
N.~Rawlinson and M.~Sambridge.
\newblock The fast marching method: An effective tool for tomographic imaging
  and tracking multiple phases in complex layered media.
\newblock {\em Exploration Geophysics}, 36:341--350, 2005.

\bibitem{ref:S96}
J.~A. Sethian.
\newblock A fast marching level set method for monotonically advancing fronts.
\newblock {\em Proceedings of the National Academy of Sciences}, 93:1591--1595,
  1996.

\bibitem{ref:S99book}
J.~A. Sethian.
\newblock {\em Level set methods and fast marching methods, evolving interfaces
  in computational geometry, fluid mechanics, computer vision, and materical
  science}.
\newblock Cambridge University Press, New York, 1999.

\bibitem{ref:S86}
H.~M. Soner.
\newblock Optimal control with state-space constraint. {II}.
\newblock {\em SIAM Journal on Control and Optimization}, 24:1110--1122, 1986.

\bibitem{ref:SA94}
P.~Spalart and S.~Allmaras.
\newblock A one-equation turbulence model for aerodynamic flows.
\newblock American Institute of Aeronautics and Astronautics 30th Aerospace
  Sciences Meeting and Exhibit, 1:5-21, 1994.

\bibitem{ref:SB13}
P.~Strachota and M.~Bene{\v{s}}.
\newblock Design and verification of the mpfa scheme for three-dimensional
  phase field model of dendritic crystal growth.
\newblock pages 459--467, Berlin, Heidelberg, 2013. Numerical Mathematics and
  Advanced Applications 2011, Springer Berlin Heidelberg.

\bibitem{ref:SLSE17}
D.~Suckart, D.~Linse, E.~Schutting, and H.~Eichlseder.
\newblock Experimental and simulative investigation of flame-wall interactions
  and quenching in spark-ignition engines.
\newblock {\em Automotive and Engine Technology}, 2(1):25--38, 2017.

\bibitem{ref:THP02}
K.~A. Tomlinson, P.~J. Hunter, and A.~J. Pullan.
\newblock A finite element method for an eikonal equation model of myocardial
  excitation wavefront propagation.
\newblock {\em SIAM Journal on Applied Mathematics}, 63:324--350, 2002.

\bibitem{ref:TCOZ03}
Y.-H.~R. Tsai, L.-T. Cheng, S.~Osher, and H.-K. Zhao.
\newblock Fast sweeping algorithms for a class of {H}amilton-{J}acobi
  equations.
\newblock {\em SIAM Journal on Numerical Analysis}, 41:673--694, 2003.

\bibitem{ref:T98}
P.~G. Tucker.
\newblock Assessment of geometric multilevel convergence robustness and a wall
  distance method for flows with multiple internal boundaries.
\newblock {\em Applied Mathematical Modelling}, 22:293--311, 1998.

\bibitem{ref:T03}
P.~G. Tucker.
\newblock Differential equation-based wall distance computation for {DES} and
  {RANS}.
\newblock {\em Journal of Computational Physics}, 190:229--248, 2003.

\bibitem{ref:T11}
P.~G. Tucker.
\newblock Hybrid {H}amilton-{J}acobi-{P}oisson wall distance function model.
\newblock {\em Computers \& Fluids}, 44:130--142, 2011.

\bibitem{ref:TRSBB05}
P.~G. Tucker, C.~L. Rumsey, P.~R. Spalart, R.~E. Bartels, and R.~T. Biedron.
\newblock Computations of wall distances based on differential equations.
\newblock {\em American Institute of Aeronautics and Astronautics Journal},
  43:539--549, 2005.

\bibitem{ref:V67}
S.~R.~S. Varadhan.
\newblock On the behavior of the fundamental solution of the heat equation with
  variable coefficients.
\newblock {\em Communications on Pure and Applied Mathematics}, 20:431--455,
  1967.

\bibitem{ref:XT10}
H.~Xia and P.~G. Tucker.
\newblock Finite volume distance field and its application to medial axis
  transforms.
\newblock {\em International Journal for Numerical Methods in Engineering},
  82:114--134, 2010.

\bibitem{ref:XT11}
H.~Xia and P.~G. Tucker.
\newblock Fast equal and biased distance fields for medial axis transform with
  meshing in mind.
\newblock {\em Applied Mathematical Modelling}, 35:5804--5819, 2011.

\bibitem{ref:Z05}
H.-K. Zhao.
\newblock Fast sweeping method for eikonal equations.
\newblock {\em Mathematics of Computation}, 74:603--627, 2005.

\bibitem{ref:Z07}
H.-K. Zhao.
\newblock Parallel implementations of the fast sweep method.
\newblock {\em Journal of Computational Mathematics}, 25:421--429, 2007.

\end{thebibliography}

\end{document}